\documentclass[12pt]{article}
\usepackage{latexsym}
\usepackage{enumerate}
\usepackage{epsf,epsfig,amsfonts,a4wide}
\usepackage{mathtools}
\usepackage{amsfonts}
\usepackage{url}
\usepackage{amsmath,amssymb,amsthm}
\usepackage{epsf,epsfig,amsfonts,graphicx}
\usepackage{a4wide}

\usepackage[centerlast]{caption}
\usepackage{graphics,amsmath,amssymb,amscd}
\usepackage{amsbsy}
\usepackage{amsthm}
\usepackage{gensymb}
\usepackage{float}
\usepackage{graphicx}
\usepackage[colorlinks=true,linkcolor=blue,citecolor=red]{hyperref}
\usepackage{multicol}

\usepackage{epstopdf}

\newcommand{\be}{\begin{equation}}
\newcommand{\ee}{\end{equation}}

\newcommand{\R}{\mathbb{R}}

\newcommand{\const}{\mathrm{const}}

\usepackage{amsthm}
\usepackage{mathdots}

\usepackage{amscd}

\newcommand{\mN}{\mathbb{N}}

\newcommand{\mR}{\mathbb{R}}

\newcommand{\mZ}{\mathbb{Z}}

\newcommand{\HH}{{\cal H}}

\newcommand{\RR}{{\cal R}}

\newcommand{\eps}{\varepsilon}
\newcommand{\ph}{\varphi}

\newcommand{\me}{\mathrm{e}}

\newcommand{\dif}{\mathrm{d}}

\newtheorem{thm}{Theorem}[section]
\newtheorem{lem}{Lemma}[section]

\newtheorem{cor}{Corollary}[section]

\title{On asymptotic description of passage through \\ a resonance in quasi-linear Hamiltonian systems}
\author{Anatoly Neishtadt$^{1,2}$\footnote{E-mail: {\tt A.Neishtadt@lboro.ac.uk}}, Tan Su$^1$\footnote{E-mail: {\tt T.Su@lboro.ac.uk}}\\
$^1$ Loughborough University, Loughborough, LE11 3TU, UK\\
$^2$ Space Research Institute, Moscow, 113997, Russia}

\begin{document}
\maketitle

\begin{abstract}
We consider a quasi-linear Hamiltonian system with one and a half degrees of freedom. The Hamiltonian of this system differs by a small, $\sim\eps$, perturbing term from the Hamiltonian of a linear oscillatory system. We consider passage through a resonance: the frequency of the latter system slowly changes with time and passes through 0. The speed of this passage is of order of $\eps$. We provide asymptotic formulas that describe effects of passage through a resonance with an accuracy $O(\eps^{\frac32})$. This is an improvement of known results by Chirikov (1959), Kevorkian (1971, 1974) and Bosley (1996). The problem under consideration is a model problem that describes passage through an isolated resonance in multi-frequency quasi-linear Hamiltonian systems. 
\end{abstract}

\vskip 20pt

\section{Introduction}

Study of passage through an isolated resonance in a multi-frequency quasi-linear Hamiltonian system can be reduced to the case of one-frequency system (see, e.g., \cite{2}). The corresponding Hamiltonian has the form
\begin{equation}\label{Ham}
H(I,\ph,\tau)=\omega(\tau)I+\eps H_1(I,\ph,\tau).
\end{equation}

Here $I$, $\ph$ mod $2\pi$ are conjugate canonical variables, $\tau$ is a slow time, $\dot\tau=\eps$, and $\eps$ is a small parameter, $0<\eps\ll1$. Equations of motion are
\begin{equation}\label{motion}
\dot I=-\eps\frac{\partial H_1(I,\ph,\tau)}{\partial\ph},\quad \dot\ph=\omega(\tau)+\eps\frac{\partial H_1(I,\ph,\tau)}{\partial I}.
\end{equation}

For $\eps=0$ we get an unperturbed system with the Hamiltonian $H_0(I,\tau)=\omega(\tau)I$ and action-angle variables $I$, $\ph$. The function $\omega$ is the frequency of the unperturbed motion. For some value of the slow time $\tau_*$, where there is a resonance, $\omega$ vanishes: $\omega(\tau_*)=0$. We assume that the resonance is non-degenerate: $\omega'_*=\omega'(\tau_*)\ne0$. Here ``prime'' denotes the derivative with respect to $\tau$. Let, for definiteness, $\omega'_*>0$. We assume that $\tau_*$ is the only resonant moment of the slow time: $\omega(\tau)$ is different from 0 at $\tau\ne\tau_*$.

Action $I$ is an adiabatic invariant: its changes along trajectory of (\ref{motion}) are small over long time intervals. For motion far from the resonance value $I$ oscillates with an amplitude $\sim\eps$. Passage through a narrow neighbourhood of the resonance leads to a change in $I$ of order $\sqrt\eps$ (so called jump of the adiabatic invariant). There is an asymptotic formula for this jump (\cite{2}, \cite{3}). Let $I_-$ and $I_+$ be values of $I$ along a trajectory of (\ref{motion}) at moments of slow time $\tau_-$ and $\tau_+$, where $\tau_-<\tau_*<\tau_+$. Then
\begin{equation}\label{difference}
I_+-I_-=-\sqrt\eps\int\limits_{-\infty}^{+\infty}\frac{\partial H_1(I_-,\ph_*+\frac{\omega'_*}{2}\theta^2,\tau_*)}{\partial\ph}\,\dif\theta+O(\eps).
\end{equation}

Here $\ph_*$ is the value of $\ph$ on the considered trajectory at $\tau=\tau_*$. There are formulas for change of the angle (phase) $\ph$ due to passage through the resonance as well \cite{2}.

One can replace in the left hand side of (\ref{difference}) $I_\pm$ with values of the improved adiabatic invariant $J_\pm$, but the error estimate in (\ref{difference}) still will be $O(\eps)$. It was suggested in \cite{1} to eliminate an asymmetry in (\ref{difference}) by replacing in the right hand side $I_-$ with $I_*$, where $I_*$ is the value of $I$ on the considered trajectory at $\tau=\tau_*$. A numerical simulation in \cite{1} shows that this symmetrization indeed improves the accuracy of formula (\ref{difference}) for $I_\pm$ replaced with $J_\pm$ considerably. It is conjectured in \cite{1} on the basis of the numerical simulation that the error term in the modified formula is $O(\eps^{\frac32})$.

In the current paper we prove this conjecture by means of a Hamiltonian adiabatic perturbation theory. We show that this improvement of accuracy occurs due to cancellations of many terms in formulas of the perturbation theory considered up to terms of 4th order in $\eps$. We obtain also formulas which describe change of phase $\ph$ due to passage through resonance with the same accuracy $O(\eps^{\frac32})$. As a result, we obtain formulas which allow to predict motion in post-resonance region with accuracy $O(\eps^{\frac32})$, provided that the motion in the pre-resonance region is known. In the last section we provide a numerical verification of these formulas.

\section{Main theorems}

We consider Hamiltonian system (\ref{motion}) with Hamilton's function (\ref{Ham}). We assume that the function $H$ is of class $C^4$ for $(I,\ph,\tau)\in D=D_I\times \R \times D_\tau$, where $D_I$ and $D_\tau$ are some open intervals in $\R$. We assume that $H_1$ is $2\pi$--periodic in $\ph$ and that the frequency $\omega$ does not vanish in $D_\tau$ other than at $\tau=\tau_*$. At the resonance state, $\omega(\tau_*)=0$, but $\omega'(\tau_*)=\omega'_*\ne 0$.

Let $I(t)$, $\ph(t)$ be a solution of (\ref{motion}) on a time interval $[t_-, t_+]$, where $t_{\pm}= \tau_{\pm}/\eps$, and $\tau_{\pm}$ are some constants. Let $\tau_-<\tau_*<\tau_+$. We denote $I_\pm=I(t_\pm)$, $\ph_\pm=\ph(t_\pm)$, and $I_*=I(t_*)$, $\ph_*=\ph(t_*)$. Here $t_*=\tau_*/\eps$.

\begin{thm}\label{thm1}
\begin{eqnarray}
I_++\eps u_1(I_+,\ph_+,\tau_+)&=&I_-+\eps u_1(I_-,\ph_-,\tau_-) -\sqrt\eps\int\limits_{-\infty}^{+\infty}\frac{\partial H_1(I_*,\ph_*+\frac{\omega_*'}{2}\theta^2,\tau_*)}{\partial\ph}\,\dif\theta + O(\eps^{\frac32}), \label{action1}\\
\ph_++\eps v_1(I_+,\ph_+,\tau_+)&=&\ph_-+\eps v_1(I_-,\ph_-,\tau_-) +\frac1\eps\int\limits_{\tau_-}^{\tau_+}\omega(\tau)\,\dif\tau +\int\limits_{\tau_-}^{\tau_*}\frac{\partial\bar H_1(J_-,\tau)}{\partial I}\,\dif\tau +\int\limits_{\tau_*}^{\tau_+}\frac{\partial\bar H_1(J_+,\tau)}{\partial I}\,\dif\tau \nonumber\\
&&{}+\sqrt\eps\int\limits_{-\infty}^{+\infty}\frac{\partial\widetilde H_1(I_*,\ph_*+\frac{\omega'_*}{2}\theta^2,\tau_*)}{\partial I}\,\dif\theta +{\rm p.v.}\,\eps\int\limits_{\tau_-}^{\tau_+}\frac{\partial{\cal R}_2(I_*,\tau)}{\partial I}\,\dif\tau \nonumber\\
&&{}-\frac{\eps^{\frac32}\ln\eps}{4\omega'_*}\frac{\partial^3\big\langle\widetilde H_1^2\big\rangle^\ph(I_*,\tau_*)}{\partial I^3} \int\limits_{-\infty}^{+\infty}\frac{\partial H_1(I_*,\ph_*+\frac{\omega'_*}{2}\theta^2,\tau_*)}{\partial\ph}\,\dif\theta+O(\eps^{\frac32}). \label{angle1}
\end{eqnarray}
\end{thm}

\begin{thm}\label{thm2}
\begin{eqnarray}
I_*&=&I_\pm\pm\frac12\sqrt\eps\int\limits_{-\infty}^{+\infty}\frac{\partial H_1(I_*,\ph_*+\frac{\omega'_*}{2}\theta^2,\tau_*)}{\partial\ph}\,\dif\theta+O(\eps), \label{action2}\\
\ph_*&=&\ph_\pm+\frac1\eps\int\limits_{\tau_\pm}^{\tau_*}\omega(\tau)\,\dif\tau +\int\limits_{\tau_\pm}^{\tau_*}\frac{\partial\bar H_1(J_\pm,\tau)}{\partial I}\,\dif\tau \mp\frac12\sqrt\eps\int\limits_{-\infty}^{+\infty}\frac{\partial\widetilde H_1(I_*,\ph_*+\frac{\omega_*'}{2}\theta^2,\tau_*)}{\partial I}\,\dif\theta \nonumber\\
&&{}+\frac{\eps\ln\eps}{4\omega'_*}\frac{\partial^2\big\langle\widetilde H_1^2\big\rangle^\ph(I_*,\tau_*)}{\partial I^2}+O(\eps). \label{angle2}
\end{eqnarray}
\end{thm}

In the above theorems, 
\begin{eqnarray}
u_1=\frac{\widetilde H_1}{\omega}, \quad \widetilde H_1=H_1-\bar H_1, \quad\bar H_1=\big\langle H_1\big\rangle^\ph, \quad v_1=-\int\limits_0^\ph\frac{\partial u_1(I,\psi,\tau)}{\partial I}\,\dif\psi +\left\langle\int\limits_0^\ph\frac{\partial u_1(I,\psi,\tau)}{\partial I}\,\dif\psi \right\rangle^\ph \nonumber, 
\end{eqnarray}
and angular brackets denote averaging with respect to $\ph$: $\big\langle f \big\rangle^\ph=\frac{1}{2\pi}\int\limits_0^{2\pi}f(\ph)\,\dif\ph$.

Moreover, $$J_\pm=I_\pm+\eps {u_1}_\pm,\quad {\cal R}_2(J_\pm,\tau)=-\frac1{2\omega(\tau)}\frac{\partial^2\big\langle\widetilde H_1^2\big\rangle^\ph(J_\pm,\tau)}{\partial I^2},\quad {u_1}_\pm=u_1(I_\pm,\ph_\pm,\tau_\pm).$$

The new results here are:
\begin{itemize}
\item estimate $O(\eps^{\frac32})$ in the first formula of Theorem \ref{thm1},
\item the last two terms and the estimate $O(\eps^{\frac32})$ in the second formula in Theorem \ref{thm1},
\item the last term and the estimate $O(\eps)$ in the second formula in Theorem \ref{thm2}.
\end{itemize}

Combining results of Theorems \ref{thm1}, \ref{thm2}, we can obtain prediction of motion in the post-resonance region with accuracy $O(\eps^{\frac32})$ as follows:

\begin{cor}\label{cor1}
\begin{eqnarray}
I_+&=&I_-+\eps u_1(I_-,\ph_-,\tau_-)-\eps u_1(I_-,\check\ph_+,\tau_+)-\sqrt\eps\int\limits_{-\infty}^{+\infty}\frac{\partial H_1(\check I_*,\check\ph_*+\frac{\omega_*'}{2}\theta^2,\tau_*)}{\partial\ph}\,\dif\theta + O(\eps^{\frac32}), \label{action3}\\
\ph_+&=&\ph_-+\eps v_1(I_-,\ph_-,\tau_-)-\eps v_1(I_-,\check\ph_+,\tau_+)+\frac1\eps\int\limits_{\tau_-}^{\tau_+}\omega(\tau)\,\dif\tau +\int\limits_{\tau_-}^{\tau_*}\frac{\partial\bar H_1(J_-,\tau)}{\partial I}\,\dif\tau +\int\limits_{\tau_*}^{\tau_+}\frac{\partial\bar H_1(\check J_+,\tau)}{\partial I}\,\dif\tau \nonumber\\
&&{}+\sqrt\eps\int\limits_{-\infty}^{+\infty}\frac{\partial\widetilde H_1(\check I_*,\check\ph_*+\frac{\omega'_*}{2}\theta^2,\tau_*)}{\partial I}\,\dif\theta +{\rm p.v.}\,\eps\int\limits_{\tau_-}^{\tau_+}\frac{\partial{\cal R}_2(I_-,\tau)}{\partial I}\,\dif\tau \nonumber\\
&&{}-\frac{\eps^{\frac32}\ln\eps}{4\omega'_*}\frac{\partial^3\big\langle\widetilde H_1^2\big\rangle^\ph(I_-,\tau_*)}{\partial I^3} \int\limits_{-\infty}^{+\infty}\frac{\partial H_1(I_-,\check\ph_*+\frac{\omega'_*}{2}\theta^2,\tau_*)}{\partial\ph}\,\dif\theta+O(\eps^{\frac32}), \label{angle3}
\end{eqnarray}
{\rm where}
\begin{eqnarray*}
\check\ph_+&=&\ph_-+\frac1\eps\int\limits_{\tau_-}^{\tau_+}\omega(\tau)\,\dif\tau +\int\limits_{\tau_-}^{\tau_+}\frac{\partial\bar H_1(I_-,\tau)}{\partial I}\,\dif\tau,\\
\check I_*&=&I_--\frac12\sqrt\eps\int\limits_{-\infty}^{+\infty}\frac{\partial H_1(I_-,\check\ph_*+\frac{\omega'_*}{2}\theta^2,\tau_*)}{\partial\ph}\,\dif\theta, \\
\check\ph_*&=&\ph_-+\frac1\eps\int\limits_{\tau_-}^{\tau_*}\omega(\tau)\,\dif\tau +\int\limits_{\tau_-}^{\tau_*}\frac{\partial\bar H_1(I_-,\tau)}{\partial I}\,\dif\tau +\frac12\sqrt\eps\int\limits_{-\infty}^{+\infty}\frac{\partial\widetilde H_1(I_-,\check{\check\ph}_*+\frac{\omega_*'}{2}\theta^2,\tau_*)}{\partial I}\,\dif\theta \\
&&{}+\frac{\eps\ln\eps}{4\omega'_*}\frac{\partial^2\big\langle\widetilde H_1^2\big\rangle^\ph(I_-,\tau_*)}{\partial I^2},\\
\check{\check\ph}_*&=&\ph_-+\frac1\eps\int\limits_{\tau_-}^{\tau_*}\omega(\tau)\,\dif\tau +\int\limits_{\tau_-}^{\tau_*}\frac{\partial\bar H_1(I_-,\tau)}{\partial I}\,\dif\tau,\\
\check J_+&=&J_--\sqrt\eps\int\limits_{-\infty}^{+\infty}\frac{\partial H_1(\check I_*,\check\ph_*+\frac{\omega'_*}{2}\theta^2,\tau_*)}{\partial\ph}\,\dif\theta.
\end{eqnarray*}
\end{cor}

\section{Procedure of adiabatic perturbation theory}

In order to get the above estimates, four steps of adiabatic perturbation theory are performed in both original Hamiltonian system (in subsection \ref{S3.1}) and an approximate Hamiltonian system (in subsection \ref{S3.2}). We follow an approach of \cite{4} here.

\subsection{Original Hamiltonian system}
\label{S3.1}

Consider dynamics described by the Hamiltonian $$H(I,\ph,\tau)=\omega(\tau)I+\eps H_1(I,\ph,\tau).$$
The frequency $\omega(\tau)$ can be expanded near the resonance as
$$\omega(\tau)=\omega'_*(\tau-\tau_*)+\frac12\omega''_*(\tau-\tau_*)^2 +O(\tau-\tau_*)^3$$
with $\omega'_*\ne0$.

So we can get formula for $\ph$ as
\begin{equation}\label{phi}
\ph=\ph_*+\frac1\eps\left(\frac12\omega'_*(\tau-\tau_*)^2 +\frac16\omega''_*(\tau-\tau_*)^3+O(\tau-\tau_*)^4\right) +\int\limits_{\tau_*}^\tau\frac{\partial H_1(I,\ph,\tau_1)}{\partial I}\,\dif\tau_1.
\end{equation}
Let
\begin{equation}\label{defS}
\eps S(J,\ph,\tau,\eps)=\eps S_1(J,\ph,\tau) +\eps^2S_2(J,\ph,\tau) +\eps^3S_3(J,\ph,\tau) +\eps^4S_4(J,\ph,\tau),
\end{equation}
where $S$ is $2\pi$-periodic in $\ph$ and $\langle S\rangle^\ph=0$. Here we have not defined $S_j$ yet.

Make the canonical transformation of variables $(I,\ph)\mapsto(J,\psi)$ with the generating function $J\ph+\eps S(J,\ph,\tau,\eps)$. Old and new variables are related as follows:
\begin{equation}\label{J}
I=J+\eps\dfrac{\partial S(J,\ph,\tau,\eps)}{\partial\ph}, \quad
\psi=\ph+\eps\dfrac{\partial S(J,\ph,\tau,\eps)}{\partial J}.
\end{equation}
The new Hamiltonian, which describes dynamics of variables $J$, $\psi$, is
\begin{eqnarray} 
\label{Hnew1}
\HH(J,\ph(J,\psi,\tau),\tau)&=&\omega(\tau)\left(J+\eps\frac{\partial S}{\partial\ph}\right) +\eps H_1(J+\eps\frac{\partial S}{\partial\ph},\ph,\tau) +\eps\frac{\partial S}{\partial t} \nonumber\\
&=&\omega(\tau)J+\eps\omega(\tau)\frac{\partial S}{\partial\ph} +\eps H_1(J,\ph,\tau) +\eps^2\frac{\partial H_1(J,\ph,\tau)}{\partial I}\frac{\partial S}{\partial\ph} +\eps^2\frac{\partial S}{\partial\tau} \nonumber\\
&&{}+\frac{\eps^3}2\frac{\partial^2H_1(J,\ph,\tau)}{\partial I^2}\left(\frac{\partial S}{\partial\ph}\right)^2 +\frac{\eps^4}{3!}\frac{\partial^3H_1(J,\ph,\tau)}{\partial I^3}\left(\frac{\partial S}{\partial\ph}\right)^3\nonumber \\
&&{}+\frac{\eps^5}{4!}\frac{\partial^4H_1(J+\theta\eps\frac{\partial S}{\partial\ph},\ph,\tau)}{\partial I^4}\left(\frac{\partial S}{\partial\ph}\right)^4,\quad 0<\theta<1.
\end{eqnarray}

We would like to find $S_j, \, j=1,..., 4$ such that there is no dependence on the new phase $\psi$ in the new Hamiltonian $\HH$ in terms up to 4th order in $\eps$. Thus the new Hamiltonian should have a form 
\begin{equation}
\label{Hnew2}
\HH(J,\ph(J,\psi,\tau),\tau)=\omega(\tau)J+\eps\RR_1(J,\tau) +\eps^2\RR_2(J,\tau) +\eps^3\RR_3(J,\tau) +\eps^4\RR_4(J,\tau) +\eps^5\RR_5(J,\psi,\tau).
\end{equation}
Here we have not defined functions $\RR_j$ yet.

Equating terms of the same order in $\eps$ in (\ref{Hnew1}) and (\ref{Hnew2})
we find that functions $\dfrac{\partial S_j}{\partial\ph}$ have forms:

\begin{equation}\label{dSdphi}
\left.\begin{split}
&\frac{\partial S_1}{\partial\ph}=\frac{\alpha^1(J,\ph,\tau)}{\omega(\tau)}, \\
&\frac{\partial S_2}{\partial\ph}=\frac{\omega'(\tau)\alpha^2(J,\ph,\tau)}{\omega^3(\tau)}+\frac{\beta^2(J,\ph,\tau)}{\omega^2(\tau)}, \\
&\frac{\partial S_3}{\partial\ph}=\frac{\big(\omega'(\tau)\big)^2\alpha^3(J,\ph,\tau)}{\omega^5(\tau)}+\frac{\beta^3(J,\ph,\tau)}{\omega^4(\tau)}, \\
&\frac{\partial S_4}{\partial\ph}=\frac{\big(\omega'(\tau)\big)^3\alpha^4(J,\ph,\tau)}{\omega^7(\tau)}+\frac{\beta^4(J,\ph,\tau)}{\omega^6(\tau)},
\end{split}\qquad\right\}
\end{equation}

where $\alpha^1$, $\alpha^2$, $\alpha^3$, $\alpha^4$, $\beta^2$, $\beta^3$, $\beta^4$ are smooth functions.

Denote $$\hat f(J,\ph,\tau)=\int f(J,\ph,\tau)\,\dif\ph-\left\langle\int f(J,\ph,\tau)\,\dif\ph\right\rangle^\ph.$$ We can find the explicit form of $S_j$:
\begin{equation}\label{S}
\left.\begin{split}
&S_1=-\frac1\omega\hat{\widetilde H}_1,\\
&S_2=-\frac{\omega'}{\omega^3}\hat{\hat{\widetilde H}}_1+\frac{\hat\beta^2}{\omega^2},\\
&S_3=-\frac{3\omega'^2}{\omega^5}\hat{\hat{\hat{\widetilde H}}}_1+\frac{\hat\beta^3}{\omega^4},\\
&S_4=-\frac{15\omega'^3}{\omega^7}\hat{\hat{\hat{\hat{\widetilde H}}}}_1 +\frac{\hat\beta^4}{\omega^6}
\end{split}\qquad\right\}
\end{equation}
satisfying $\langle S_j\rangle^\ph=0$, $j=1,2,3,4$.

So we have the expressions
\begin{equation}\label{dSdJ}
\left.\begin{split}
\eps\frac{\partial S_1}{\partial J}&=-\frac\eps\omega\frac{\partial\hat{\widetilde H}_1}{\partial J},\\
\eps^2\frac{\partial S_2}{\partial J}&=-\eps^2\frac{\omega'}{\omega^3}\frac{\partial\hat{\hat{\widetilde H}}_1}{\partial J}+\frac{\eps^2}{\omega^2}\frac{\partial\hat\beta^2}{\partial J},\\
\eps^3\frac{\partial S_3}{\partial J}&=-\eps^3\frac{3\omega'^2}{\omega^5}\frac{\partial\hat{\hat{\hat{\widetilde H}}}_1}{\partial J}+\frac{\eps^3}{\omega^4}\frac{\partial\hat\beta^3}{\partial J},\\
\eps^4\frac{\partial S_4}{\partial J}&=-\eps^4\frac{15\omega'^3}{\omega^7}\frac{\partial\hat{\hat{\hat{\hat{\widetilde H}}}}_1}{\partial J} +\frac{\eps^4}{\omega^6}\frac{\partial\hat\beta^4}{\partial J}.
\end{split}\qquad\right\}
\end{equation}

Also
\begin{eqnarray*}
&&\RR_1(J,\tau)=\big\langle H_1(J,\ph,\tau)\big\rangle^\ph, \\
&&\RR_2(J,\tau)=\left\langle\frac{\partial H_1(J,\ph,\tau)}{\partial I} \frac{\partial S_1}{\partial\ph}\right\rangle^\ph, \\
&&\RR_3(J,\tau)=\left\langle\frac12\frac{\partial^2H_1(J,\ph,\tau)}{\partial I^2} \left(\frac{\partial S_1}{\partial\ph}\right)^2+\frac{\partial H_1(J,\ph,\tau)}{\partial I}\frac{\partial S_2}{\partial\ph}\right\rangle^\ph, \\
&&\RR_4(J,\tau)=\left\langle\frac16\frac{\partial^3H_1(J,\ph,\tau)}{\partial I^3} \left(\frac{\partial S_1}{\partial\ph}\right)^3 +\frac{\partial^2H_1(J,\ph,\tau)}{\partial I^2}\frac{\partial S_1}{\partial\ph}\frac{\partial S_2}{\partial\ph} +\frac{\partial H_1(J,\ph,\tau)}{\partial I}\frac{\partial S_3}{\partial\ph}\right\rangle^\ph.
\end{eqnarray*}

The new Hamiltonian is
$$\HH(J,\ph(J,\psi,\tau),\tau)=\omega(\tau)J+\eps\RR_1+\eps^2\RR_2+\eps^3\RR_3+\eps^4\RR_4 +\frac{\eps^5O(1)}{\omega^7(\tau)}+\eps^5\frac{\partial S_4}{\partial\tau}.$$
For $|\tau-\tau_*|\ge\sqrt\eps$, the motion is described by differential equations:
\begin{eqnarray}\label{dotJ}
&&\dot J=-\dfrac{\partial\HH}{\partial\ph}\cdot\dfrac{\partial\ph}{\partial\psi}=-\dfrac{\partial\HH}{\partial\ph}\cdot\left(\dfrac{\partial\psi}{\partial\ph}\right)^{-1}\nonumber\\
&&\quad=-\dfrac{\partial\HH}{\partial\ph}\left(1+\eps\dfrac{\gamma^1(J,\ph,\tau)}{\omega(\tau)}+\eps^2\dfrac{\gamma^2(J,\ph,\tau)}{\omega^3(\tau)}+\eps^3\dfrac{\gamma^3(J,\ph,\tau)}{\omega^5(\tau)}+\eps^4\dfrac{\gamma^4(J,\ph,\tau)}{\omega^7(\tau)}\right)^{-1} \nonumber\\
&&\quad=\left(\eps^5\dfrac{7\big(\omega'(\tau)\big)^4\alpha^4(J,\ph,\tau)}{\omega^8(\tau)} +\eps^5\dfrac{\gamma(J,\ph,\tau)}{\omega^7(\tau)}\right) \left(1-\displaystyle\sum\limits_{k=1}^{4}\eps^k\dfrac{\gamma^k(J,\ph,\tau)}{\omega^{2k-1}(\tau)}+O\left(\dfrac{\eps^2}{(\tau-\tau_*)^2}\right)\right), \nonumber\\\\
&&\dot \psi=\dfrac{\partial\HH}{\partial J}+\dfrac{\partial\HH}{\partial\ph}\cdot\dfrac{\partial\ph}{\partial J}=\dfrac{\partial\HH}{\partial J}+\dfrac{\partial\HH}{\partial\ph}\cdot\left(-\eps\dfrac{\partial^2S}{\partial J^2}\right)\left(1+\eps\dfrac{\partial^2S}{\partial J\partial\ph}\right)^{-1}\nonumber\\
&&\quad=\dfrac{\partial\HH}{\partial J}-\dfrac{\partial\HH}{\partial\ph}\cdot\dfrac{\partial}{\partial J}\left(\eps\dfrac{\partial S_1}{\partial J}+\eps^2\dfrac{\partial S_2}{\partial J}+\eps^3\dfrac{\partial S_3}{\partial J}+\eps^4\dfrac{\partial S_4}{\partial J}\right)\left(1+O\bigg(\dfrac\eps\omega\bigg)\right)^{-1}\nonumber\\
&&\quad=\omega(\tau)+\dfrac{\partial}{\partial J}\left(\eps\RR_1+\eps^2\RR_2+\eps^3\RR_3+\eps^4\RR_4 +\dfrac{\eps^5O(1)}{\omega^7(\tau)}+\eps^5\dfrac{\partial S_4}{\partial\tau}\right)\nonumber\\
&&\quad\quad{}+\left(\eps^5\dfrac{7\omega'^4\alpha^4}{\omega^8}+\eps^5\dfrac{\gamma}{\omega^7}\right)\left(\eps\dfrac{\widetilde\gamma^1}{\omega}+\eps^2\dfrac{\widetilde\gamma^2}{\omega^3}+\eps^3\dfrac{\widetilde\gamma^3}{\omega^5}+\eps^4\dfrac{\widetilde\gamma^4}{\omega^7}\right)\left(1+O\bigg(\dfrac\eps\omega\bigg)\right)\nonumber\\
&&\quad=\omega(\tau)+\dfrac{\partial}{\partial J}\left(\eps\RR_1+\eps^2\RR_2+\eps^3\RR_3+\eps^4\RR_4+\eps^5\dfrac{\partial S_4}{\partial\tau}\right)+\displaystyle\sum\limits_{k=5}^{9}\eps^k\dfrac{\gamma^k(J,\ph,\tau)}{\omega^{2k-3}(\tau)}+O\left(\dfrac{\eps^7}{(\tau-\tau_*)^{10}}\right).\nonumber\\
\end{eqnarray}
Here $\gamma$, $\gamma^1,\ldots,\gamma^9$, $\widetilde\gamma^1,\ldots,\widetilde\gamma^4$ are smooth functions.

{\bf Remark.} By differentiating $\psi=\ph+\eps\dfrac{\partial S(J,\ph(J,\psi,\tau),\tau)}{\partial J}$ with respect to $J$ on both sides, we obtain $0=\dfrac{\partial\ph}{\partial J}+\eps\dfrac{\partial^2S}{\partial J^2}+\eps\dfrac{\partial^2S}{\partial J\partial\ph}\dfrac{\partial\ph}{\partial J}$. Therefore, $\dfrac{\partial\ph}{\partial J}=-\eps\dfrac{\partial^2S}{\partial J^2}\left(1+\dfrac{\partial^2S}{\partial J\partial\ph}\right)^{-1}$.

\subsection{Approximate Hamiltonian system}
\label{S3.2}
Now let us consider the approximate Hamiltonian
$$H(I_a,\ph_a,\tau)=\Omega(\tau)I_a+\eps H_1(I_*,\ph_a,\tau_*)$$
Equations of motion are
$$ 
\dot I_a=-\eps\dfrac{\partial H_1(I_*,\ph_a,\tau_*)}{\partial\ph}, \quad \dot\ph_a=\Omega(\tau).
$$
Here $\Omega(\tau)=\omega'_*(\tau-\tau_*)$, $\omega'_*\ne0$.

We will consider the solution of these equations with initial conditions at resonance, i.e. when $\tau=\tau_*$: $I_a(\tau_*)=I_*=I(t_*)$, $\ph_a(\tau_*)=\ph_*=\ph(t_*)$. We get the formula for $\ph_a$ as
\begin{equation}\label{phia}
\ph_a=\ph_*+\frac1{2\eps}\omega'_*(\tau-\tau_*)^2.
\end{equation}

Let
\begin{equation}\label{defSa}
\eps S^a(\ph_a,\tau,\eps)=\eps S^a_1(\ph_a,\tau)+\eps^2S^a_2(\ph_a,\tau) +\eps^3S^a_3(\ph_a,\tau)+\eps^4S^a_4(\ph_a,\tau)
\end{equation}
where $S^a$ is $2\pi$-periodic in $\ph_a$ and $\big\langle S^a\big\rangle^{\ph_a}=0$. Here we have not defined $S^a_j$ yet.

Make the canonical transformation of variables $(I_a,\ph_a)\mapsto(J_a,\psi_a)$ with a generating function $J_a\ph_a+\eps S^a(\ph_a,\tau,\eps)$. The old and new variables are related as follows:
\begin{equation}\label{Ja}
I_a=J_a+\eps\dfrac{\partial S^a(\ph_a,\tau,\eps)}{\partial\ph_a}, \quad
\psi_a=\ph_a.
\end{equation}

The new Hamiltonian, which describes dynamics of variables $J_a$, $\ph_a$, is
\begin{eqnarray}
\label{HHnew1}
\HH^a(J_a,\ph_a,\tau)&=&\Omega(\tau)\left(J_a+\eps\frac{\partial S^a}{\partial\ph_a}\right)+\eps H_1(I_*,\ph_a,\tau_*)+\eps\frac{\partial S^a}{\partial t} \nonumber\\
&=&\Omega(\tau)J_a+\eps\Omega(\tau)\frac{\partial S^a}{\partial\ph_a}+\eps H_1(I_*,\ph_a,\tau_*)+\eps^2\frac{\partial S^a}{\partial\tau}.
\end{eqnarray}
We would like to find $S^a_j, \, j=1,..., 4$ such that there is no dependence on the phase $\ph_a$ in the new Hamiltonian $\HH^a$ in terms up to 4th order in $\eps$. Thus the new Hamiltonian should have a form 
\begin{equation}
\label{HHnew2}
\HH^a(J_a,\ph_a,\tau) =\Omega(\tau)J_a+\eps\RR_1^a(J_a,\tau)+\eps^2\RR_2^a(J_a,\tau)+\eps^3\RR_3^a(J_a,\tau)+\eps^4\RR_4^a(J_a,\tau)+\eps^5\RR_5^a(J_a,\ph_a,\tau).
\end{equation}

Equating terms of the same order in $\eps$ in (\ref{HHnew1}) and (\ref{HHnew2})
we find that functions $\dfrac{\partial S^a_j}{\partial\ph_a}$ have forms
\begin{equation}\label{dSadphi}
\left.\begin{split}
&\frac{\partial S_1^a}{\partial\ph_a}=\frac{\alpha^1(I_*,\ph_a,\tau_*)}{\Omega(\tau)}, \\
&\frac{\partial S_2^a}{\partial\ph_a}=\frac{\Omega'(\tau)\alpha^2(I_*,\ph_a,\tau_*)}{\Omega^3(\tau)}, \\
&\frac{\partial S_3^a}{\partial\ph_a}=\frac{\big(\Omega'(\tau)\big)^2\alpha^3(I_*,\ph_a,\tau_*)}{\Omega^5(\tau)}, \\
&\frac{\partial S_4^a}{\partial\ph_a}=\frac{\big(\Omega'(\tau)\big)^3\alpha^4(I_*,\ph_a,\tau_*)}{\Omega^7(\tau)}.
\end{split}\qquad\right\}
\end{equation}

The new Hamiltonian is
$$\HH^a(J_a,\ph_a,\tau)=\Omega(\tau)J_a+\eps\big\langle H_1(I_*,\ph_a,\tau_*)\big\rangle^{\ph_a} +\eps^5\frac{\partial S^a_4(\ph_a,\tau)}{\partial\tau}.$$
The motion is described by differential equations
\begin{equation}
\begin{cases}
\ \dot J_a=-\eps^5\dfrac{\partial^2 S^a_4(\ph_a,\tau)}{\partial\ph_a\partial\tau} =\eps^5\dfrac{7\big(\Omega'(\tau)\big)^4\alpha^4(I_*,\ph_a,\tau_*)}{\Omega^8(\tau)} +\eps^5\dfrac{\widetilde\gamma(I_*,\ph_a,\tau_*)}{\Omega^7(\tau)}, \\\\
\ \dot \ph_a=\Omega(\tau).
\end{cases}
\end{equation}
Here $\widetilde\gamma$ is a smooth function.

\section{Proofs of the theorems}

We will prove asymptotic formulas for the action variable $I$ in both Theorems \ref{thm1} and \ref{thm2} first, and then asymptotic formulas for the angle variable $\ph$ in these Theorems. Denote $\tau_l=\tau_* - \eps^{1/2}$, $\tau_r=\tau_* + \eps^{1/2}$, $t_{l,r}
=\tau_{l,r}/ \eps^{1/2}$. Denote $u=-\frac{\partial S}{\partial\ph}$ and $u^a=-\frac{\partial S^a}{\partial\ph_a}$. In the following text, we use the notations $Z_j=Z(t_j)$, where $Z=I,I_a, J, J_a, \ph, \ph_a, \psi, S, S^a, S_k, S_k^a, u, u^a, u_k, u_k^a, v, v_k$ and $t_j=t_-, t_l, t_*, t_r, t_+$, $k=1,2,3,4$.

For simplicity of the exposition we will assume that $\tau_+$ and $\tau_-$ are symmetric with respect to $\tau_*$: $\tau_+-\tau_*=\tau_*-\tau_-$. General case can be easily reduced to this one.

\subsection{Principal lemmas}
The following lemmas will be used in the proof.

\begin{lem}\label{err_omega}
$$\frac{\big(\omega'(\tau)\big)^a}{\omega^b(\tau)}=\frac{\big(\Omega'(\tau)\big)^a}{\Omega^b(\tau)}+\frac{O(1)}{(\tau-\tau_*)^{b-1}}=\frac{\const}{(\tau-\tau_*)^b}+\frac{O(1)}{(\tau-\tau_*)^{b-1}},$$
where $a,b\in\mathbb{N}$.
\end{lem}

\begin{lem}\label{est}
\begin{eqnarray*}
&&I_l-I_*=O(\sqrt\eps),\quad I_r-I_*=O(\sqrt\eps),\quad I_l-I_r=O(\sqrt\eps).\\
&&J_l-I_*=O(\sqrt\eps),\quad J_r-I_*=O(\sqrt\eps),\quad J_l-J_r=O(\sqrt\eps).\\
&&\ph_l-\ph_{a_l}=O(\sqrt\eps),\quad \ph_r-\ph_{a_r}=O(\sqrt\eps),\quad \ph_l-\ph_r=O(\sqrt\eps).
\end{eqnarray*}
\end{lem}

\begin{lem}\label{est_comb}
$$(J_l-I_*)+(J_r-I_*)=O(\eps),\quad (\ph_l-\ph_{a_l}) +(\ph_r-\ph_{a_r})=O(\eps).$$
\end{lem}

\proof

\quad$\bullet$ For $\big[(J_l-I_*)+(J_r-I_*)\big]$, substitute with formulas (\ref{J}) and (\ref{dSdphi}):
\begin{eqnarray*}
&&\quad(J_l-I_*)+(J_r-I_*)=\left(I_l-I_*-\eps\frac{\partial S(J_l,\ph_l,\tau_l)}{\partial\ph}\right)+\left(I_r-I_*-\eps\frac{\partial S(J_r,\ph_r,\tau_r)}{\partial\ph}\right)\\
&&=\eps\int\limits_{t_l}^{t_*}\frac{\partial H_1(I,\ph,\tau)}{\partial\ph}\,\dif t -\eps\frac{\hat\alpha^1(J_l,\ph_l,\tau_l)}{\tau_l-\tau_*}-\eps^2\frac{\hat\alpha^2(J_l,\ph_l,\tau_l)}{(\tau_l-\tau_*)^3}-\eps^3\frac{\hat\alpha^3(J_l,\ph_l,\tau_l)}{(\tau_l-\tau_*)^5}-\eps^4\frac{\hat\alpha^4(J_l,\ph_l,\tau_l)}{(\tau_l-\tau_*)^7}+O(\eps)\\
&&{}-\eps\int\limits_{t_*}^{t_r}\frac{\partial H_1(I,\ph,\tau)}{\partial\ph}\,\dif t -\eps\frac{\hat\alpha^1(J_r,\ph_r,\tau_r)}{\tau_r-\tau_*}-\eps^2\frac{\hat\alpha^2(J_r,\ph_r,\tau_r)}{(\tau_r-\tau_*)^3}-\eps^3\frac{\hat\alpha^3(J_r,\ph_r,\tau_r)}{(\tau_r-\tau_*)^5}-\eps^4\frac{\hat\alpha^4(J_r,\ph_r,\tau_r)}{(\tau_r-\tau_*)^7}+O(\eps)\\
&&=\int\limits_{\tau_l}^{\tau_*}\frac{\partial H_1(I_*,\ph_a,\tau_*)}{\partial\ph}\,\dif\tau-\int\limits_{\tau_*}^{\tau_r}\frac{\partial H_1(I_*,\ph_a,\tau_*)}{\partial\ph}\,\dif\tau-\sum_{k=1}^{4}\eps^k\frac{\hat\alpha^k(J_l,\ph_l,\tau_l)}{(\tau_l-\tau_*)^{2k-1}}-\sum_{k=1}^{4}\eps^k\frac{\hat\alpha^k(J_l,\ph_l,\tau_l)}{(\tau_r-\tau_*)^{2k-1}}+O(\eps)\\
&&=O(\eps).
\end{eqnarray*}

\quad$\bullet$ For $\big[(\ph_l-\ph_{a_l}) +(\ph_r-\ph_{a_r})\big]$, substitute with formulas (\ref{phi}) and (\ref{phia}):
\begin{eqnarray*}
&&\quad(\ph_l-\ph_{a_l}) +(\ph_r-\ph_{a_r})\\
&&=\frac{\omega''_*}{6\eps}\big[(\tau_l-\tau_*)^3 +(\tau_r-\tau_*)^3\big]+\frac1\eps O(\tau_l-\tau_*)^4+\int\limits_{\tau_*}^{\tau_l}\frac{\partial H_1(I,\ph,\tau_1)}{\partial I}\,\dif\tau_1+\int\limits_{\tau_*}^{\tau_r}\frac{\partial H_1(I,\ph,\tau_1)}{\partial I}\,\dif\tau_1\\
&&=\int\limits_{\tau_*}^{\tau_l}\frac{\partial H_1(I_*,\ph_a,\tau_*)}{\partial I}\,\dif\tau_1+\int\limits_{\tau_*}^{\tau_r}\frac{\partial H_1(I_*,\ph_a,\tau_*)}{\partial I}\,\dif\tau_1+O(\eps)\\
&&=O(\eps)
\end{eqnarray*}
\qed

\begin{lem}\label{est_frac}
$$\eps^a\frac{(\tau-\tau_*)^b}{\omega^c(\tau)}f_0(I,\ph,\tau)=O(\eps^{a+\frac{b-c}2}),$$
if $\tau\in\{\tau_l,\tau_r\}$. Here $a\in\mR$, $b,c\in\mZ$, and $f_0$ is a smooth function. 

The same estimate is valid if $\omega(\tau)$ is replaced with $\Omega(\tau)$.
\end{lem}

\begin{lem}\label{est_int}
$$\int\limits_{t_l}^t\eps^a\frac{f_0(I(t_1),\ph(t_1),\eps t_1)}{\omega^c(\eps t_1)}\,\dif t_1=O(\eps^{a-\frac{c+1}2}),$$
if $t\in[t_-,t_l]$. Here $a\in\mR$, $c\in\mZ$, $c\ge2$, and $f_0$ is a smooth function. The same estimate is valid if $t_l$ is replaced with $t_r$ and $t\in[t_r,t_+]$.

The same estimate is valid if $\omega(\tau)$ is replaced with $\Omega(\tau)$.
\end{lem}

\begin{lem}[Cancellation lemma near the resonance]
\label{cancNEAR}
$$\int\limits_{\tau_l}^{\tau_r}(\tau-\tau_*)^{2k-1}f_1(I_*,\ph_a,\tau_*)\,\dif\tau=0$$
where $k\in\mathbb{N}$ and $f_1$ is a smooth function.
\end{lem}

\begin{lem}[Cancellation lemma far from the resonance on symmetric intervals]
\label{cancFAR}
$$\int\limits_{\tau_-}^{\tau_l}(\tau-\tau_*)^{2k-1}f_2(I_*,\ph_a,\tau_*)\,\dif\tau+\int\limits_{\tau_r}^{\tau_+}(\tau-\tau_*)^{2k-1}f_2(I_*,\ph_a,\tau_*)\,\dif\tau=0$$
where $k\in\mathbb{Z}$ and $f_2$ is a smooth function.
\end{lem}

\begin{lem}\label{varbI} Let $f=f(I,\ph,\tau)$ be a twice continuously differentiable function. Then
\begin{eqnarray*}
&&\quad f(I,\ph,\tau)-f(I_*,\ph_a,\tau_*)\\
&&=-\frac{\partial f(I_*,\ph_a,\tau_*)}{\partial I}\int\limits_{\tau_*}^{\tau}\frac{\partial H_1(I_*,\ph_a,\tau_*)}{\partial\ph}\,\dif\tau_1 +\frac{\partial f(I_*,\ph_a,\tau_*)}{\partial \ph}\left(\frac{\omega''_*(\tau-\tau_*)^3}{6\eps}+\int\limits_{\tau_*}^{\tau}\frac{\partial H_1(I_*,\ph_a,\tau_*)}{\partial I}\,\dif\tau_1\right)\\
&&\quad{}+\frac{\partial f(I_*,\ph_a,\tau_*)}{\partial \tau}(\tau-\tau_*)+O(\tau-\tau_*)^2+O\left(\frac{(\tau-\tau_*)^4}{\eps}\right)+O\left(\frac{(\tau-\tau_*)^6}{\eps^2}\right)\,.
\end{eqnarray*}
\end{lem}

\proof
\begin{eqnarray*}
&&\quad f(I,\ph,\tau)-f(I_*,\ph_a,\tau_*)\\
&&=\frac{\partial f(I_*,\ph_a,\tau_*)}{\partial I}(I-I_*)+\frac{\partial f(I_*,\ph_a,\tau_*)}{\partial \ph}(\ph-\ph_a)+\frac{\partial f(I_*,\ph_a,\tau_*)}{\partial \tau}(\tau-\tau_*)\\
&&\quad{}+O(I-I_*)^2+O(\ph-\ph_a)^2+O(\tau-\tau_*)^2\\
&&=-\frac{\partial f(I_*,\ph_a,\tau_*)}{\partial I}\int\limits_{\tau_*}^{\tau}\frac{\partial H_1(I,\ph,\tau_1)}{\partial\ph}\,\dif\tau_1 +\frac{\partial f(I_*,\ph_a,\tau_*)}{\partial \ph}\left(\frac{\omega''_*(\tau-\tau_*)^3}{6\eps}+\int\limits_{\tau_*}^{\tau}\frac{\partial H_1(I,\ph,\tau_1)}{\partial I}\,\dif\tau_1\right)\\
&&\quad{}+\frac{\partial f(I_*,\ph_a,\tau_*)}{\partial \tau}(\tau-\tau_*)+O(\tau-\tau_*)^2+O\left(\frac{(\tau-\tau_*)^4}{\eps}\right)+O\left(\frac{(\tau-\tau_*)^6}{\eps^2}\right)\\
&&=-\frac{\partial f(I_*,\ph_a,\tau_*)}{\partial I}\int\limits_{\tau_*}^{\tau}\frac{\partial H_1(I_*,\ph_a,\tau_*)}{\partial\ph}\,\dif\tau_1 +\frac{\partial f(I_*,\ph_a,\tau_*)}{\partial \ph}\left(\frac{\omega''_*(\tau-\tau_*)^3}{6\eps}+\int\limits_{\tau_*}^{\tau}\frac{\partial H_1(I_*,\ph_a,\tau_*)}{\partial I}\,\dif\tau_1\right)\\
&&\quad{}+\frac{\partial f(I_*,\ph_a,\tau_*)}{\partial \tau}(\tau-\tau_*)+O(\tau-\tau_*)^2+O\left(\frac{(\tau-\tau_*)^4}{\eps}\right)+O\left(\frac{(\tau-\tau_*)^6}{\eps^2}\right)
\end{eqnarray*}
\qed

\begin{lem}\label{varbJ} Let $f=f(J,\ph,\tau)$ be a twice continuously differentiable function. Then
for $\tau\in[\tau_-,\tau_l]\cup[\tau_r,\tau_+]$,
\begin{eqnarray*}
&&\quad f(J,\ph,\tau)-f(I_*,\ph_a,\tau_*)\\
&&=-\frac{\partial f(I_*,\ph_a,\tau_*)}{\partial I}\left(\,\int\limits_{\tau_*}^{\tau}\frac{\partial H_1(I_*,\ph_a,\tau_*)}{\partial\ph}\,\dif\tau_1 \right) \\
&&\quad{}+\frac{\partial f(I_*,\ph_a,\tau_*)}{\partial \ph}\left(\frac{\omega''_*(\tau-\tau_*)^3}{6\eps}+\int\limits_{\tau_*}^{\tau}\frac{\partial H_1(I_*,\ph_a,\tau_*)}{\partial I}\,\dif\tau_1\right)+\frac{\partial f(I_*,\ph_a,\tau_*)}{\partial \tau}(\tau-\tau_*)\\
&&\quad{}+O(\eps)+O(\tau-\tau_*)^2+O\left(\frac{(\tau-\tau_*)^4}{\eps}\right)+O\left(\frac{(\tau-\tau_*)^6}{\eps^2}\right)\\
&&\quad{}+\eps\frac{f^1(I_*,\ph_a,\tau_*)}{\tau-\tau_*}+\eps^2\frac{f^2(I_*,\ph_a,\tau_*)}{(\tau-\tau_*)^3}+\eps^3\frac{f^3(I_*,\ph_a,\tau_*)}{(\tau-\tau_*)^5}+\eps^4\frac{f^4(I_*,\ph_a,\tau_*)}{(\tau-\tau_*)^7}
\end{eqnarray*}
where $f^1$, $f^2$, $f^3$, $f^4$ are smooth functions.
\end{lem}

\proof
\begin{eqnarray*}
&&\quad f(J,\ph,\tau)-f(I_*,\ph_a,\tau_*)\\
&&=\frac{\partial f(I_*,\ph_a,\tau_*)}{\partial I}(J-I_*)+\frac{\partial f(I_*,\ph_a,\tau_*)}{\partial \ph}(\ph-\ph_a)+\frac{\partial f(I_*,\ph_a,\tau_*)}{\partial \tau}(\tau-\tau_*)\\
&&\quad{}+O(J-I_*)^2+O(\ph-\ph_a)^2+O(\tau-\tau_*)^2\\
&&=\frac{\partial f(I_*,\ph_a,\tau_*)}{\partial I}\left(I-I_*-\eps\frac{\partial S(J,\ph,\tau)}{\partial\ph}\right)+\frac{\partial f(I_*,\ph_a,\tau_*)}{\partial \ph}(\ph-\ph_a)+\frac{\partial f(I_*,\ph_a,\tau_*)}{\partial \tau}(\tau-\tau_*)\\
&&\quad{}+O(J-I_*)^2+O(\ph-\ph_a)^2+O(\tau-\tau_*)^2\\
&&=\frac{\partial f(I_*,\ph_a,\tau_*)}{\partial I}\left(I-I_*-\eps\frac{\partial S(J,\ph,\tau)}{\partial\ph}\right)+\frac{\partial f(I_*,\ph_a,\tau_*)}{\partial \ph}
\left(\frac{\omega''_*(\tau-\tau_*)^3}{6\eps}+\int\limits_{\tau_*}^{\tau}\frac{\partial H_1(I_*,\ph_a,\tau_*)}{\partial I}\,\dif\tau_1\right)\\
&&\quad{}+ \frac{\partial f(I_*,\ph_a,\tau_*)}{\partial \tau}(\tau-\tau_*)+
O(\eps)+O(\tau-\tau_*)^2+O\left(\frac{(\tau-\tau_*)^4}{\eps}\right)+O\left(\frac{(\tau-\tau_*)^6}{\eps^2}\right).
\end{eqnarray*}
Then
\begin{eqnarray*}
&&I-I_*=-\int\limits_{\tau_*}^{\tau}\frac{\partial H_1(I,\ph,\tau)}{\partial\ph}\,\dif\tau_1 =-\int\limits_{\tau_*}^{\tau}\frac{\partial H_1(I_*,\ph_a,\tau_*)}{\partial\ph}\,\dif\tau_1+O(\tau-\tau_*)^2+O\left(\frac{(\tau-\tau_*)^4}{\eps}\right).
\end{eqnarray*}
Also
\begin{eqnarray*}
&&\quad\eps\frac{\partial S(J,\ph,\tau)}{\partial\ph}= \eps\frac{\partial S(I_*,\ph_a,\tau)}{\partial\ph}+O(\eps)+O(\tau-\tau_*)^2\\
&&=-\eps\frac{f^1(I_*,\ph_a,\tau_*)}{\tau-\tau_*}-\eps^2\frac{f^2(I_*,\ph_a,\tau_*)}{(\tau-\tau_*)^3}-\eps^3\frac{f^3(I_*,\ph_a,\tau_*)}{(\tau-\tau_*)^5}-\eps^4\frac{f^4(I_*,\ph_a,\tau_*)}{(\tau-\tau_*)^7}+O(\eps)+O(\tau-\tau_*)^2.
\end{eqnarray*}
This implies the result of the lemma.
\qed

\subsection{Proof of formulas for action variable}

\subsubsection{Principal identity for formula (\ref{action1})}
Let
\begin{eqnarray*}
e_1=\int\limits_{t_r}^{t_+}(\dot J-\dot J_a)\,\dif t,&\quad &e_2=\int\limits_{t_l}^{t_r}(\dot I-\dot I_a)\,\dif t, \\
e_3=\int\limits_{t_-}^{t_l}(\dot J-\dot J_a)\,\dif t,&\quad &e_4=-\int\limits_{-\infty}^{t_-}\dot I_a\,\dif t-\int\limits_{t_+}^{+\infty}\dot I_a\,\dif t+\eps(u^a_+-u^a_-).
\end{eqnarray*}
We have the identity:
\begin{eqnarray*}
I_++\eps u_+&=&(J_+)\ =\ J_r+\int\limits_{t_r}^{t_+}\dot J\,\dif t\\
&=&I_r+\eps u_r+\int\limits_{t_r}^{t_+}(\dot J-\dot J_a)\,\dif t+\int\limits_{t_r}^{t_+}\dot J_a\,\dif t \\
&=&I_r+\eps u_r+e_1+\int\limits_{t_r}^{t_+}(\dot I_a+\eps\dot u^a)\,\dif t \\
&=&I_r+\eps u_r+e_1+\int\limits_{t_r}^{t_+}\dot I_a\,\dif t+\eps u^a_+-\eps u^a_r \\
&=&I_l+\int\limits_{t_l}^{t_r}\dot I\,\dif t+\int\limits_{t_r}^{t_+}\dot I_a\,\dif t+\eps u_r+\eps u^a_+-\eps u^a_r+e_1 \\
&=&I_l+\eps u_l+\int\limits_{t_l}^{t_r}(\dot I-\dot I_a)\,\dif t+\int\limits_{t_l}^{t_r}\dot I_a\,\dif t+\int\limits_{t_r}^{t_+}\dot I_a\,\dif t+\eps u_r-\eps u_l+\eps u^a_+-\eps u^a_r+e_1 \\
&=&I_-+\eps u_-+\int\limits_{t_-}^{t_l}\dot J\,\dif t+\int\limits_{t_l}^{t_r}\dot I_a\,\dif t+\int\limits_{t_r}^{t_+}\dot I_a\,\dif t+\eps u_r-\eps u_l+\eps u^a_+-\eps u^a_r+e_1+e_2 \\
&=&I_-+\eps u_-+\int\limits_{t_-}^{t_l}\dot J_a\,\dif t+\int\limits_{t_l}^{t_r}\dot I_a\,\dif t+\int\limits_{t_r}^{t_+}\dot I_a\,\dif t+\eps u_r-\eps u_l+\eps u^a_+-\eps u^a_r+e_1+e_2+e_3 \\
&=&I_-+\eps u_-+\int\limits_{t_-}^{t_+}\dot I_a\,\dif t+\eps u^a_+-\eps u^a_-+\eps u_r-\eps u_l+\eps u^a_l-\eps u^a_r+e_1+e_2+e_3 \\
&=&I_-+\eps u_-+\int\limits_{-\infty}^{+\infty}\dot I_a\,\dif t+\eps u_r-\eps u_l+\eps u^a_l-\eps u^a_r+e_1+e_2+e_3+e_4.
\end{eqnarray*}

We should prove that $\eps u_r-\eps u_l+\eps u^a_l-\eps u^a_r=O(\eps^{\frac32})$, $e_1+e_3= O(\eps^{\frac32})$, and $e_j= O(\eps^{\frac32}),\, j=2,4$.

\subsubsection{Estimate of $\eps u_r-\eps u_l+\eps u^a_l-\eps u^a_r$}

Here we use the sum $u=u_1+u_2+u_3+u_4$, $u^a=u^a_1+u^a_2+u^a_3+u^a_4$, $u_k=-\frac{\partial S_k}{\partial\ph}$, $u^a_k=-\frac{\partial S^a_k}{\partial\ph_a}$, $k=1,2,3,4$. For convenience, we consider $S_1(J,\ph,\tau)$ as the main part of expression of $S(J,\ph,\tau)$, then consider the others.
\begin{eqnarray*}
&&\eps {u_1}_r-\eps {u_1}_l+\eps {u^a_1}_l-\eps {u^a_1}_r\\
&=&-\eps\frac{\partial {S_1}_r}{\partial\ph}+\eps\frac{\partial {S_1}_l}{\partial\ph} -\eps\frac{\partial {S^a_1}_l}{\partial\ph_a}+\eps\frac{\partial {S^a_1}_r}{\partial\ph_a} \\
&=&-\eps\frac{\alpha^1(J_r,\ph_r,\tau_r)}{\omega(\tau_r)}+\eps\frac{\alpha^1(J_l,\ph_l,\tau_l)}{\omega(\tau_l)}-\eps\frac{\alpha^1(I_*,\ph_{a_l},\tau_*)}{\Omega(\tau_l)}+\eps\frac{\alpha^1(I_*,\ph_{a_r},\tau_*)}{\Omega(\tau_r)} \\
&=&\eps\frac{\Omega(\tau_l)\alpha^1(J_l,\ph_l,\tau_l)-\omega(\tau_l)\alpha^1(I_*,\ph_{a_l},\tau_*)}{\omega(\tau_l)\Omega(\tau_l)}-\eps\frac{\Omega(\tau_r)\alpha^1(J_r,\ph_r,\tau_r)-\omega(\tau_r)\alpha^1(I_*,\ph_{a_r},\tau_*)}{\omega(\tau_r)\Omega(\tau_r)}\\
&=&\eps\frac{\big[\omega(\tau_l)-\frac12\omega''_*(\tau_l-\tau_*)^2 +O(\tau_l-\tau_*)^3\big]\alpha^1(J_l,\ph_l,\tau_l) -\omega(\tau_l)\alpha^1(I_*,\ph_{a_l},\tau_*)}{\omega(\tau_l)\Omega(\tau_l)}\\
&&-\eps\frac{\big[\omega(\tau_r)-\frac12\omega''_*(\tau_r-\tau_*)^2 +O(\tau_r-\tau_*)^3\big]\alpha^1(J_r,\ph_r,\tau_r)-\omega(\tau_r)\alpha^1(I_*,\ph_{a_r},\tau_*)}{\omega(\tau_r)\Omega(\tau_r)}\\
&=&\eps\frac{\omega(\tau_l)\big[\alpha^1(J_l,\ph_l,\tau_l)-\alpha^1(I_*,\ph_{a_l},\tau_*)\big]-\frac{\omega''_*(\tau_l-\tau_*)^2}2\alpha^1(J_l,\ph_l,\tau_l)+O(\tau_l-\tau_*)^3}{\omega(\tau_l)\Omega(\tau_l)}\\
&&-\eps\frac{\omega(\tau_r)\big[\alpha^1(J_r,\ph_r,\tau_r)-\alpha^1(I_*,\ph_{a_r},\tau_*)\big]-\frac{\omega''_*(\tau_r-\tau_*)^2}2\alpha^1(J_r,\ph_r,\tau_r)+O(\tau_r-\tau_*)^3}{\omega(\tau_r)\Omega(\tau_r)}\\
&=&\eps\frac{\alpha^1(J_l,\ph_l,\tau_l)-\alpha^1(I_*,\ph_{a_l},\tau_*)}{\Omega(\tau_l)} -\eps\frac{\omega''_*(\tau_l-\tau_*)^2\alpha^1(J_l,\ph_l,\tau_l)}{2\omega(\tau_l)\Omega(\tau_l)}\\
&&-\eps\frac{\alpha^1(J_r,\ph_r,\tau_r)-\alpha^1(I_*,\ph_{a_r},\tau_*)}{\Omega(\tau_r)} +\eps\frac{\omega''_*(\tau_r-\tau_*)^2\alpha^1(J_r,\ph_r,\tau_r)}{2\omega(\tau_r)\Omega(\tau_r)}+O(\eps^{\frac32})\\
&=&\frac{\eps}{\Omega(\tau_l)}\Big[\big(\alpha^1(J_l,\ph_l,\tau_l)-\alpha^1(I_*,\ph_{a_l},\tau_*)\big)+\big(\alpha^1(J_r,\ph_r,\tau_r)-\alpha^1(I_*,\ph_{a_r},\tau_*)\big)\Big]\\
&&-\frac\eps2\frac{\omega''_*(\tau_l-\tau_*)^2}{\Omega(\tau_l)}\left(\frac{\alpha^1(J_l,\ph_l,\tau_l)}{\omega(\tau_l)}+\frac{\alpha^1(J_r,\ph_r,\tau_r)}{\omega(\tau_r)}\right)+O(\eps^{\frac32})\\
&=&O(\sqrt\eps)E_{1+}+O(\eps^{\frac32})\left(\frac{\alpha^1(J_l,\ph_l,\tau_l)}{\omega(\tau_l)}+\frac{\alpha^1(J_l,\ph_l,\tau_l)+O(\sqrt\eps)}{\omega(\tau_r)}\right)+O(\eps^{\frac32})\\
&=&O(\sqrt\eps)E_{1+}+O(\eps^{\frac32})\left(\alpha^1(J_l,\ph_l,\tau_l)\cdot\frac{\omega(\tau_l)+\omega(\tau_r)}{\omega(\tau_l)\omega(\tau_r)}\right)+O(\eps^{\frac32})\\
&=&O(\sqrt\eps)E_{1+}+O(\eps^{\frac32})
\end{eqnarray*}

Here $E_{1+}=\big[\alpha^1(J_l,\ph_l,\tau_l)-\alpha^1(I_*,\ph_{a_l},\tau_*)\big]+\big[\alpha^1(J_r,\ph_r,\tau_r)-\alpha^1(I_*,\ph_{a_r},\tau_*)\big]$, and we are going to estimate $E_{1+}$ using Lemmas \ref {est} and \ref{est_comb}:

\begin{eqnarray*}
E_{1+}&=&\frac{\partial\alpha^1(I_*,\ph_{a_l},\tau_*)}{\partial I}(J_l-I_*) +\frac{\partial\alpha^1(I_*,\ph_{a_l},\tau_*)}{\partial \ph}(\ph_l-\ph_{a_l}) +\frac{\partial\alpha^1(I_*,\ph_{a_l},\tau_*)}{\partial \tau}(\tau_l-\tau_*)\\
&&{}+\frac{\partial\alpha^1(I_*,\ph_{a_r},\tau_*)}{\partial I}(J_r-I_*) +\frac{\partial\alpha^1(I_*,\ph_{a_r},\tau_*)}{\partial \ph}(\ph_r-\ph_{a_r}) +\frac{\partial\alpha^1(I_*,\ph_{a_r},\tau_*)}{\partial \tau}(\tau_r-\tau_*)\\
&&{}+O(J_l-I_*)^2+O(\ph_l-\ph_{a_l})^2+O(\tau_l-\tau_*)^2+O(J_r-I_*)^2+O(\ph_r-\ph_{a_r})^2+O(\tau_r-\tau_*)^2\\
&=&\frac{\partial\alpha^1(I_*,\ph_{a_l},\tau_*)}{\partial I}\big[(J_l-I_*)+(J_r-I_*)\big]+\frac{\partial\alpha^1(I_*,\ph_{a_l},\tau_*)}{\partial \ph}\big[(\ph_l-\ph_{a_l})+(\ph_r-\ph_{a_r})\big]\\
&&{}+\frac{\partial\alpha^1(I_*,\ph_{a_l},\tau_*)}{\partial \tau}\big[(\tau_l-\tau_*)+(\tau_r-\tau_*)\big]+O(\ph_{a_l}-\ph_{a_r})^2+O(\eps)\\
&=&O(\eps)
\end{eqnarray*}
Therefore, $\eps {u_1}_r-\eps {u_1}_l+\eps {u^a_1}_l-\eps {u^a_1}_r=O(\eps^{\frac32})$.

Then we consider the term $S_2(J,\ph,\tau)$ in expression of $S(J,\ph,\tau)$. Here 
\begin{eqnarray*}
&&\frac{\partial S_2}{\partial\ph}=\frac{\omega'(\tau)\alpha^2(J,\ph,\tau)}{\omega^3(\tau)}+\frac{\beta^2(J,\ph,\tau)}{\omega^2(\tau)}, \\
&&\frac{\partial S_2^a}{\partial\ph_a}=\frac{\Omega'(\tau)\alpha^2(I_*,\ph_a,\tau_*)}{\Omega^3(\tau)}.
\end{eqnarray*}
\begin{eqnarray*}
&&\eps {u_2}_r-\eps {u_2}_l+\eps {u^a_2}_l-\eps {u^a_2}_r\\
&=&-\eps^2\frac{\partial {S_2}_r}{\partial\ph}+\eps^2\frac{\partial {S_2}_l}{\partial\ph} -\eps^2\frac{\partial {S^a_2}_l}{\partial\ph_a}+\eps^2\frac{\partial {S^a_2}_r}{\partial\ph_a} \\
&=&-\eps^2\frac{\omega'(\tau_r)\alpha^2(J_r,\ph_r,\tau_r)}{\omega^3(\tau_r)}-\eps^2\frac{\beta^2(J_r,\ph_r,\tau_r)}{\omega^2(\tau_r)}+\eps^2\frac{\omega'(\tau_l)\alpha^2(J_l,\ph_l,\tau_l)}{\omega^3(\tau_l)}+\eps^2\frac{\beta^2(J_l,\ph_l,\tau_l)}{\omega^2(\tau_l)}\\
&&{}-\eps^2\frac{\Omega'(\tau_l)\alpha^2(I_*,\ph_{a_l},\tau_*)}{\Omega^3(\tau_l)}+\eps^2\frac{\Omega'(\tau_r)\alpha^2(I_*,\ph_{a_r},\tau_*)}{\Omega^3(\tau_r)} \\
\end{eqnarray*}

From Lemmas \ref{est} and \ref{est_frac},
\begin{eqnarray*}
&&-\eps^2\frac{\beta^2(J_r,\ph_r,\tau_r)}{\omega^2(\tau_r)}+\eps^2\frac{\beta^2(J_l,\ph_l,\tau_l)}{\omega^2(\tau_l)}\\
&=&-\eps^2\frac{\beta^2(J_r,\ph_r,\tau_r)}{\Omega^2(\tau_r)}+\eps^2\frac{\beta^2(J_l,\ph_l,\tau_l)}{\Omega^2(\tau_l)}+O(\eps^{\frac32})\\
&=&\frac{\eps^2}{\Omega^2(\tau_l)}\big[\beta^2(J_l,\ph_l,\tau_l)-\beta^2(J_r,\ph_r,\tau_r)\big]+O(\eps^{\frac32})\\
&=&\frac{\eps^2}{\Omega^2(\tau_l)}\left[\frac{\partial\beta^2}{\partial J}(J_l-J_r)+\frac{\partial\beta^2}{\partial\ph}(\ph_l-\ph_r)+\frac{\partial\beta^2}{\partial\tau}(\tau_l-\tau_r)\right]+O(\eps^{\frac32})\\
&=&O(\eps^{\frac32})
\end{eqnarray*}

Also we know that
\begin{eqnarray*}
\Omega^3(\tau)&=&\big[\omega(\tau)-\frac12\omega''_*(\tau-\tau_*)^2 +O(\tau-\tau_*)^3\big]^3\\
&=&\omega^3(\tau)-\frac32\omega''_*\omega^2(\tau)(\tau-\tau_*)^2+O(\tau-\tau_*)^5.
\end{eqnarray*}

So
\begin{eqnarray*}
&&\eps {u_r}_2-\eps {u_l}_2+\eps {u_{a_l}}_2-\eps {u_{a_r}}_2\\
&=&\eps^2\frac{\Omega^3(\tau_l)\omega'(\tau_l)\alpha^2(J_l,\ph_l,\tau_l) -\omega^3(\tau_l)\Omega'(\tau_l)\alpha^2(I_*,\ph_{a_l},\tau_*)}{\omega^3(\tau_l)\Omega^3(\tau_l)}\\
&&{}-\eps^2\frac{\Omega^3(\tau_r)\omega'(\tau_r)\alpha^2(J_r,\ph_r,\tau_r) -\omega^3(\tau_r)\Omega'(\tau_r)\alpha^2(I_*,\ph_{a_r},\tau_*)}{\omega^3(\tau_r)\Omega^3(\tau_r)}+O(\eps^{\frac32})\\
&=&\eps^2\frac{\big[\omega^3(\tau_l)-\frac32\omega''_*\omega^2(\tau_l)(\tau_l-\tau_*)^2+O(\tau_l-\tau_*)^5\big]\omega'(\tau_l)\alpha^2(J_l,\ph_l,\tau_l)}{\omega^3(\tau_l)\Omega^3(\tau_l)}\\
&&{}-\eps^2\frac{\omega^3(\tau_l)\big[\omega'(\tau_l)-\omega''_*(\tau_l-\tau_*) +O(\tau_l-\tau_*)^2\big] \alpha^2(I_*,\ph_{a_l},\tau_*)}{\omega^3(\tau_l)\Omega^3(\tau_l)}\\
&&{}-\eps^2\frac{\big[\omega^3(\tau_r)-\frac32\omega''_*\omega^2(\tau_r)(\tau_r-\tau_*)^2+O(\tau_r-\tau_*)^5\big]\omega'(\tau_r)\alpha^2(J_r,\ph_r,\tau_r)}{\omega^3(\tau_r)\Omega^3(\tau_r)}\\
&&{}+\eps^2\frac{\omega^3(\tau_r)\big[\omega'(\tau_r)-\omega''_*(\tau_r-\tau_*) +O(\tau_r-\tau_*)^2\big] \alpha^2(I_*,\ph_{a_r},\tau_*)}{\omega^3(\tau_r)\Omega^3(\tau_r)}+O(\eps^{\frac32})\\
&=&\eps^2\frac{\omega^3(\tau_l)\omega'(\tau_l)\big[\alpha^2(J_l,\ph_l,\tau_l)-\alpha^2(I_*,\ph_{a_l},\tau_*)\big]}{\omega^3(\tau_l)\Omega^3(\tau_l)}\\
&&{}-\eps^2\frac{\frac32\omega''_*\omega'(\tau_l)\omega^2(\tau_l)(\tau_l-\tau_*)^2\alpha^2(J_l,\ph_l,\tau_l)}{\omega^3(\tau_l)\Omega^3(\tau_l)} +\eps^2\frac{\omega''_*\omega^3(\tau_l)\alpha^2(I_*,\ph_{a_l},\tau_*)(\tau_l-\tau_*)}{\omega^3(\tau_l)\Omega^3(\tau_l)}\\
&&{}-\eps^2\frac{\omega^3(\tau_r)\omega'(\tau_r)\big[\alpha^2(J_r,\ph_r,\tau_r)-\alpha^2(I_*,\ph_{a_r},\tau_*)\big]}{\omega^3(\tau_r)\Omega^3(\tau_r)}\\
&&{}+\eps^2\frac{\frac32\omega''_*\omega'(\tau_r)\omega^2(\tau_r)(\tau_r-\tau_*)^2\alpha^2(J_r,\ph_r,\tau_r)}{\omega^3(\tau_r)\Omega^3(\tau_r)} -\eps^2\frac{\omega''_*\omega^3(\tau_r)\alpha^2(I_*,\ph_{a_r},\tau_*)(\tau_r-\tau_*)}{\omega^3(\tau_r)\Omega^3(\tau_r)}+O(\eps^{\frac32})\\
&=&\eps^2\frac{\omega'_*}{\Omega^3(\tau_l)}\big[\alpha^2(J_l,\ph_l,\tau_l)-\alpha^2(I_*,\ph_{a_l},\tau_*)\big]-\eps^2\frac{\omega'_*}{\Omega^3(\tau_r)}\big[\alpha^2(J_r,\ph_r,\tau_r)-\alpha^2(I_*,\ph_{a_r},\tau_*)\big] \\
&&{}-\frac32\eps^2\frac{\omega''_*\omega'(\tau_l)}{\omega(\tau_l)\Omega^3(\tau_l)}(\tau_l-\tau_*)^2\alpha^2(J_l,\ph_l,\tau_l) +\frac32\eps^2\frac{\omega''_*\omega'(\tau_r)}{\omega(\tau_r)\Omega^3(\tau_r)}(\tau_r-\tau_*)^2\alpha^2(J_r,\ph_r,\tau_r) \\
&&{}+\eps^2\frac{\omega''_*}{\Omega^3(\tau_l)}(\tau_l-\tau_*)\alpha^2(I_*,\ph_{a_l},\tau_*)-\eps^2\frac{\omega''_*}{\Omega^3(\tau_r)}(\tau_r-\tau_*)\alpha^2(I_*,\ph_{a_r},\tau_*)+O(\eps^{\frac32})\\
\end{eqnarray*}
\begin{eqnarray*}
&=&\eps^2\frac{\omega'_*}{\Omega^3(\tau_l)}\Big[\big(\alpha^2(J_l,\ph_l,\tau_l)-\alpha^2(I_*,\ph_{a_l},\tau_*)\big)+\big(\alpha^2(J_r,\ph_r,\tau_r)-\alpha^2(I_*,\ph_{a_r},\tau_*)\big)\Big]\\
&&{}+\eps^2\frac{\omega''_*(\tau_l-\tau_*)}{\Omega^3(\tau_l)}\big[\alpha^2(J_l,\ph_l,\tau_l)-\alpha^2(I_*,\ph_{a_l},\tau_*)\big]-\eps^2\frac{\omega''_*(\tau_r-\tau_*)}{\Omega^3(\tau_r)}\big[\alpha^2(J_r,\ph_r,\tau_r)-\alpha^2(I_*,\ph_{a_r},\tau_*)\big]\\
&&{}-\frac32\eps^2\frac{\omega''_*\omega'_*(\tau_l-\tau_*)^2}{\omega(\tau_l)\Omega^3(\tau_l)}\big[\alpha^2(J_l,\ph_l,\tau_l)-\alpha^2(J_r,\ph_r,\tau_r)\big] \\
&&{}+\eps^2\frac{\omega''_*(\tau_l-\tau_*)}{\Omega^3(\tau_l)}\big[\alpha^2(I_*,\ph_{a_l},\tau_*)-\alpha^2(I_*,\ph_{a_r},\tau_*)\big] +O(\eps^{\frac32})\\
&=&O(\sqrt\eps)E_{2+}+O(\eps)\big[\alpha^2(J_l,\ph_l,\tau_l)-\alpha^2(I_*,\ph_{a_l},\tau_*)\big]+O(\eps)\big[\alpha^2(J_r,\ph_r,\tau_r)-\alpha^2(I_*,\ph_{a_r},\tau_*)\big]\\
&&{}+O(\eps)\big[\alpha^2(J_l,\ph_l,\tau_l)-\alpha^2(J_r,\ph_r,\tau_r)\big]+O(\eps)\big[\alpha^2(I_*,\ph_{a_l},\tau_*)-\alpha^2(I_*,\ph_{a_r},\tau_*)\big]+O(\eps^{\frac32})\\
&=&O(\sqrt\eps)E_{2+}+O(\eps^{\frac32}),
\end{eqnarray*}
where $E_{2+}=\big[\alpha^2(J_l,\ph_l,\tau_l)-\alpha^2(I_*,\ph_{a_l},\tau_*)\big]+\big[\alpha^2(J_r,\ph_r,\tau_r)-\alpha^2(I_*,\ph_{a_r},\tau_*)\big]$.

Similarly to $E_{1+}=O(\eps)$, we can obtain $E_{2+}=O(\eps)$ from Lemma \ref{est_comb}. Therefore, it is true that $\eps {u_2}_r-\eps {u_2}_l+\eps {u^a_2}_l-\eps {u^a_2}_r=O(\eps^{\frac32})$.

Similarly, we can derive that 
$$\eps {u_{3,4}}_r-\eps {u_{3,4}}_l+\eps {u^a_{3,4}}_l-\eps {u^a_{3,4}}_r=O(\eps^{\frac32}).$$
Therefore,
$$\eps u_r-\eps u_l+\eps u^a_l-\eps u^a_r=O(\eps^{\frac32}).$$

\subsubsection{Estimate of $e_2$}
We simply use Lemma \ref{varbI} and Lemma \ref{cancNEAR} in order to get estimate:
$$e_2=\int\limits_{t_l}^{t_r}(\dot I-\dot I_a)\,\dif t=-\eps\int\limits_{t_l}^{t_r}\left(\frac{\partial H_1(I,\ph,\tau)}{\partial\ph}-\frac{\partial H_1(I_*,\ph_a,\tau_*)}{\partial\ph}\right)\,\dif t=O(\eps^{\frac32}).$$

\subsubsection{Estimate of $e_1+e_3$}
For combined term $e_1+e_3$, we consider $\dot J=\dot J_1-\dot J_2+O\big(\frac{\eps^7}{(\tau-\tau_*)^{10}}\big)$, where 
\begin{eqnarray*}
&&\dot J_1=\eps^5\dfrac{7\big(\omega'(\tau)\big)^4\alpha^4(J,\ph,\tau)}{\omega^8(\tau)} +\eps^5\dfrac{\gamma(J,\ph,\tau)}{\omega^7(\tau)}, \\
&&\dot J_2= \left(\eps^5\dfrac{7\big(\omega'(\tau)\big)^4\alpha^4(J,\ph,\tau)}{\omega^8(\tau)} +\eps^5\dfrac{\gamma(J,\ph,\tau)}{\omega^7(\tau)}\right) \cdot \sum\limits_{k=1}^{4}\eps^k\frac{\gamma^k(J,\ph,\tau)}{\omega^{2k-1}(\tau)}.
\end{eqnarray*}
and also
$$\dot J_a=\eps^5\frac{7\big(\Omega'(\tau)\big)^4\alpha^4(I_*,\ph_a,\tau_*)}{\Omega^8(\tau)}+\eps^5\frac{\widetilde\gamma(I_*,\ph_a,\tau_*)}{\Omega^7(\tau)}.$$

For $\dot J_1-\dot J_a$, we apply Lemma \ref{err_omega}, then Lemma \ref{varbJ}:
\begin{eqnarray*}
\dot J_1-\dot J_a&=&\eps^5\frac{7\big(\Omega'(\tau)\big)^4}{\Omega^8(\tau)} \big[\alpha^4(J,\ph,\tau)-\alpha^4(I_*,\ph_a,\tau_*)\big]+\eps^5\frac{\check\gamma(J,\ph,\tau)}{(\tau-\tau_*)^7}-\eps^5\frac{\widetilde\gamma(I_*,\ph_a,\tau_*)}{\Omega^7(\tau)}+O\left(\frac{\eps^5}{(\tau-\tau_*)^6}\right)\\
&=&\frac{\eps^5\delta_1(I_*,\ph_a,\tau_*)}{(\tau-\tau_*)^8}\int\limits_{\tau_*}^{\tau}\frac{\partial H_1(I_*,\ph_a,\tau_*)}{\partial\ph}\,\dif\tau_1 \\
&&{}+\frac{\eps^4\delta_2(I_*,\ph_a,\tau_*)}{(\tau-\tau_*)^5} +\frac{\eps^5\delta_3(I_*,\ph_a,\tau_*)}{(\tau-\tau_*)^8}\int\limits_{\tau_*}^{\tau}\frac{\partial H_1(I_*,\ph_a,\tau_*)}{\partial I}\,\dif\tau_1 +\frac{\eps^5\delta_4(I_*,\ph_a,\tau_*)}{(\tau-\tau_*)^7} \\
&&{}+O\left(\frac{\eps^7}{(\tau-\tau_*)^{10}}\right) +O\left(\frac{\eps^6}{(\tau-\tau_*)^8}\right) +O\left(\frac{\eps^5}{(\tau-\tau_*)^6}\right) +O\left(\frac{\eps^4}{(\tau-\tau_*)^4}\right) +O\left(\frac{\eps^3}{(\tau-\tau_*)^2}\right)\\
&&{}+\frac{\eps^6\delta_5(I_*,\ph_a,\tau_*)}{(\tau-\tau_*)^9} +\frac{\eps^7\delta_6(I_*,\ph_a,\tau_*)}{(\tau-\tau_*)^{11}} +\frac{\eps^8\delta_7(I_*,\ph_a,\tau_*)}{(\tau-\tau_*)^{13}} +\frac{\eps^9\delta_8(I_*,\ph_a,\tau_*)}{(\tau-\tau_*)^{15}}.
\end{eqnarray*}
Here $\check\gamma$ is a smooth function.
Also
\begin{eqnarray*}
\dot J_2&=& \left(\eps^5\dfrac{7\big(\omega'(\tau)\big)^4\alpha^4(J,\ph,\tau)}{\omega^8(\tau)} +\eps^5\dfrac{\gamma(J,\ph,\tau)}{\omega^7(\tau)}\right) \cdot \sum\limits_{k=1}^{4}\eps^k\frac{\gamma^k(J,\ph,\tau)}{\omega^{2k-1}(\tau)}\\
&=&\sum\limits_{k=1}^{4}\left(\eps^{k+5}\frac{\widetilde\gamma^k(J,\ph,\tau)}{(\tau-\tau_*)^{2k+7}}+O\left(\frac{\eps^{k+5}}{(\tau-\tau_*)^{2k+6}}\right)\right)\\
&=&\sum\limits_{k=1}^{4}\left(\eps^{k+5}\frac{\widetilde\gamma^k(I_*,\ph_a,\tau_*)}{(\tau-\tau_*)^{2k+7}} +O\left(\frac{\eps^{k+5}}{(\tau-\tau_*)^{2k+6}}\right) +O\left(\frac{\eps^{k+6}}{(\tau-\tau_*)^{2k+8}}\right) +O\left(\frac{\eps^{k+4}}{(\tau-\tau_*)^{2k+4}}\right)\right)\\
&=&\sum\limits_{k=1}^{4}\eps^{k+5}\frac{\widetilde\gamma^k(I_*,\ph_a,\tau_*)}{(\tau-\tau_*)^{2k+7}} +\sum\limits_{k=1}^{6}O\left(\frac{\eps^{k+4}}{(\tau-\tau_*)^{2k+4}}\right)
\end{eqnarray*}

Using Lemma \ref{est_int} for estimate of integrals and Lemma \ref{cancFAR} for cancellation in symetric intervals, we obtain:
$$e_1+e_3=\int\limits_{t_r}^{t_+}(\dot J_1-\dot J_a)\,\dif t +\int\limits_{t_-}^{t_l}(\dot J_1-\dot J_a)\,\dif t -\int\limits_{t_r}^{t_+}\dot J_2\,\dif t -\int\limits_{t_-}^{t_l}\dot J_2\,\dif t+O(\eps^{\frac32})=O(\eps^{\frac32}).$$

\subsubsection{Estimate of $e_4$}\label{e4}

For term $e_4$, we apply an integration by parts:
\begin{eqnarray*}
-\int\limits_{-\infty}^{t_-}\dot I_a\,\dif t&=&\eps \int\limits_{-\infty}^{t_-}\frac{\partial\widetilde H_1(I_*,\ph_a,\tau_*)}{\partial\ph}\,\dif t\\
&=&\eps\int\limits_{-\infty}^{t_-}\frac{\partial\widetilde H_1(I_*,\ph_a,\tau_*)}{\partial\ph}\cdot\dot\ph_a\cdot\frac1{\Omega(\tau)}\,\dif t\\
&=&\eps\int\limits_{-\infty}^{t_-}\frac1{\Omega(\tau)}\,\dif\widetilde H_1\\
&=&\frac{\eps\widetilde H_1(I_*,\ph_a,\tau_*)}{\Omega(\eps t)}\Bigg|_{-\infty}^{t_-} -\eps\int\limits_{-\infty}^{t_-}\widetilde H_1\cdot\left(-\frac{\eps\Omega'(\tau)}{\Omega^2(\tau)}\right)\,\dif t\\
&=&\frac{\eps\widetilde H_1(I_*,{\ph_a}_-,\tau_*)}{\Omega(\tau_-)} +\eps^2\int\limits_{-\infty}^{t_-}\frac{\widetilde H_1(I_*,\ph_a,\tau_*)}{\omega'_*(\tau-\tau_*)^2}\,\dif t\\
&=&\frac{\eps\widetilde H_1(I_*,{\ph_a}_-,\tau_*)}{\Omega(\tau_-)} +O(\eps^2)\\
&=&\eps u^a_-+O(\eps^2),
\end{eqnarray*}
and similarly $$-\int\limits^{+\infty}_{t_+}\dot I_a\,\dif t=-\frac{\eps\widetilde H_1(I_*,{\ph_a}_+,\tau_*)}{\Omega(\tau_+)} +O(\eps^2)=-\eps u^a_++O(\eps^2).$$

Therefore, the estimate $e_4=O(\eps^2)$ is obtained.

By joining estimate of terms $e_j$ together, and taking into account that $u_\pm={u_1}_\pm+O(\eps)$, as well as the identity
$$\int\limits_{-\infty}^{+\infty}\dot I_a\,\dif t = -\eps\int\limits_{-\infty}^{+\infty}\frac{\partial H_1(I_*,\ph_a,\tau_*)}{\partial\ph}\,\dif t=-\sqrt\eps\int\limits_{-\infty}^{+\infty}\frac{\partial H_1(I_*,\ph_*+\frac{\omega_*'}{2}\theta^2,\tau_*)}{\partial\ph}\,\dif\theta$$
where $\theta=\dfrac{\tau-\tau_*}{\sqrt\eps}$, we have finished the proof of the first formula of Theorem \ref{thm1}.

\subsubsection{Principal identity for formula (\ref{action2})}
Let
$$\tilde e_1=\int\limits_{t_l}^{t_*}(\dot I-\dot I_a)\,\dif t,\quad \tilde e_2=\int\limits_{t_-}^{t_l}(\dot J-\dot J_a)\,\dif t,\quad \tilde e_3=-\int\limits_{-\infty}^{t_-}\dot I_a\,\dif t.$$

We have
\begin{eqnarray*}
I_*&=&I_l+\int\limits_{t_l}^{t_*}\dot I\,\dif t\\
&=&I_l+\int\limits_{t_l}^{t_*}(\dot I-\dot I_a)\,\dif t+\int\limits_{t_l}^{t_*}\dot I_a\,\dif t \\
&=&J_l-\eps u_l+\int\limits_{t_l}^{t_*}\dot I_a\,\dif t+\tilde e_1\\
&=&J_-+\int\limits_{t_-}^{t_l}\dot J\,\dif t-\eps u_l+\int\limits_{t_l}^{t_*}\dot I_a\,\dif t+\tilde e_1\\
&=&I_-+\eps u_-+\int\limits_{t_-}^{t_l}(\dot J-\dot J_a)\,\dif t +\int\limits_{t_-}^{t_l}\dot J_a\,\dif t-\eps u_l +\int\limits_{t_l}^{t_*}\dot I_a\,\dif t +\tilde e_1\\
&=&I_-+\int\limits_{t_-}^{t_l}(\dot I_a+\eps\dot u^a)\,\dif t+\eps u_--\eps u_l +\int\limits_{t_l}^{t_*}\dot I_a\,\dif t +\tilde e_1+\tilde e_2\\
&=&I_-+\int\limits_{t_-}^{t_*}\dot I_a\,\dif t+\eps u^a_l-\eps u^a_-+\eps u_--\eps u_l +\tilde e_1+\tilde e_2\\
&=&I_-+\int\limits_{-\infty}^{t_*}\dot I_a\,\dif t-\int\limits_{-\infty}^{t_-}\dot I_a\,\dif t+\eps u^a_l-\eps u^a_-+\eps u_--\eps u_l +\tilde e_1+\tilde e_2\\
&=&I_--\eps\int\limits_{-\infty}^{t_*}\frac{\partial H_1(I_*,\ph_a,\tau_*)}{\partial\ph}\,\dif t+\eps u^a_l-\eps u^a_-+\eps u_--\eps u_l +\tilde e_1+\tilde e_2+\tilde e_3\\
&=&I_--\frac12\sqrt\eps\int\limits_{-\infty}^{+\infty}\frac{\partial H_1(I_*,\ph_*+\frac{\omega'_*}{2}\theta^2,\tau_*)}{\partial\ph}\,\dif\theta+\eps u^a_l-\eps u^a_-+\eps u_--\eps u_l +\tilde e_1+\tilde e_2+\tilde e_3.
\end{eqnarray*}

Making use of Lemmas \ref{err_omega}, \ref{est}, \ref{est_frac}, \ref{est_int}, \ref{varbI}, \ref{varbJ}, one can show that $\tilde e_j=O(\eps)$, $j=1,2,3$, $\eps u^a_l-\eps u_l=O(\eps)$, $\eps u_-=O(\eps)$, $\eps u^a_-=O(\eps)$. This leads to the first formula of Theorem \ref{thm2} with the sign ``$-$". The proof of the formula with the sign ``$+$" is completely analogous.

\subsection{Proof of formulas for angle variable}

\subsubsection{Principal lemmas}

\begin{lem}\label{replaceNEAR}
For $k\in\mN$ and smooth function $f(J,\ph,\tau)$,
\begin{eqnarray*}
&(a)&\left.\frac{\eps^kf(J,\ph,\tau)}{\omega^{2k-2}(\tau)}\right|_{t_l}^{t_r}=O(\eps^{\frac32}),\\
&(b)&\left.\frac{\eps^kf(J,\ph,\tau)}{\omega^{2k-1}(\tau)}\right|_{t_l}^{t_r}=\left.\frac{\eps^kf(I_*,\ph_a,\tau_*)}{\Omega^{2k-1}(\tau)}\right|_{t_l}^{t_r}+O(\eps^{\frac32}).
\end{eqnarray*}
\end{lem}

\proof
For $\tau\in\{\tau_l,\tau_r\}$, Lemma \ref{varbJ} can be simplified as $f(J,\ph,\tau)=f(I_*,\ph_a,\tau_*)+O(\tau-\tau_*)$. Thus with Lemmas \ref{err_omega} and \ref{est_frac},
\begin{eqnarray*}
(\rm a)&\left.\dfrac{\eps^kf(J,\ph,\tau)}{\omega^{2k-2}(\tau)}\right|_{t_l}^{t_r}&=\left.\frac{\eps^kf(J,\ph,\tau)}{\Omega^{2k-2}(\tau)}\right|_{t_l}^{t_r}+\left.\frac{\eps^kO(1)}{\Omega^{2k-3}(\tau)}\right|_{t_l}^{t_r}\\
&&=\left.\frac{\eps^kf(J,\ph,\tau)}{\Omega^{2k-2}(\tau)}\right|_{t_l}^{t_r}+O(\eps^{\frac32})\\
&&=\left.\frac{\eps^kf(I_*,\ph_a,\tau_*)}{\Omega^{2k-2}(\tau)}\right|_{t_l}^{t_r}+\left.\frac{\eps^kO(1)}{\Omega^{2k-3}(\tau)}\right|_{t_l}^{t_r}+O(\eps^{\frac32})\\
&&=O(\eps^{\frac32}).\\
\end{eqnarray*}
Applying the result of (a) and Lemma \ref{varbJ}, we obtain
\begin{eqnarray*}
(\rm b)&\left.\dfrac{\eps^kf(J,\ph,\tau)}{\omega^{2k-1}(\tau)}\right|_{t_l}^{t_r}&=\left.\frac{\eps^kf(J,\ph,\tau)}{\Omega^{2k-1}(\tau)}\right|_{t_l}^{t_r}+\left.\frac{\eps^k\tilde f(J,\ph,\tau)}{\Omega^{2k-2}(\tau)}\right|_{t_l}^{t_r}+\left.\frac{\eps^kO(1)}{\Omega^{2k-3}(\tau)}\right|_{t_l}^{t_r}\\
&&=\left.\frac{\eps^kf(J,\ph,\tau)}{\Omega^{2k-1}(\tau)}\right|_{t_l}^{t_r}+O(\eps^{\frac32})\\
&&=\left.\frac{\eps^kf(I_*,\ph_a,\tau_*)}{\Omega^{2k-1}(\tau)}\right|_{t_l}^{t_r}+\left.\frac{\eps^k\tilde{\tilde f}(I_*,\ph_a,\tau_*)}{\Omega^{2k-2}(\tau)}\right|_{t_l}^{t_r}+O(\eps^{\frac32})\\
&&=\left.\frac{\eps^kf(I_*,\ph_a,\tau_*)}{\Omega^{2k-1}(\tau)}\right|_{t_l}^{t_r}+O(\eps^{\frac32}).
\end{eqnarray*}
Here $\tilde f$ and $\tilde{\tilde f}$ is smooth functions.
\qed

\begin{lem}\label{replaceFAR}
Let $f=f(J,\ph,\tau)$ be a twice continuously differentiable function. Then for $k\in\mN$, $k\ge5$,
\begin{eqnarray*}
&(a)&\int\limits_{t_-}^{t_l}\frac{\eps^kf(J,\ph,\tau)}{\omega^{2k-3}(\tau)}\,\dif t +\int\limits_{t_r}^{t_+}\frac{\eps^kf(J,\ph,\tau)}{\omega^{2k-3}(\tau)}\,\dif t =O(\eps^{\frac32}),\\
&(b)&\int\limits_{t_-}^{t_l}\frac{\eps^kf(J,\ph,\tau)}{\omega^{2k-2}(\tau)}\,\dif t +\int\limits_{t_r}^{t_+}\frac{\eps^kf(J,\ph,\tau)}{\omega^{2k-2}(\tau)}\,\dif t =\int\limits_{t_-}^{t_l}\frac{\eps^kf(I_*,\ph_a,\tau_*)}{\Omega^{2k-2}(\tau)}\,\dif t +\int\limits_{t_r}^{t_+}\frac{\eps^kf(I_*,\ph_a,\tau_*)}{\Omega^{2k-2}(\tau)}\,\dif t +O(\eps^{\frac32}).
\end{eqnarray*}
\end{lem}

\proof
For $t\in[t_-,t_l]\cup[t_r,t_+]$, we can simplify Lemma \ref{varbJ} as follows:
$$f(J,\ph,\tau)=f(I_*,\ph_a,\tau_*)+O(\tau-\tau_*)+O\left(\frac{(\tau-\tau_*)^3}{\eps}\right)+O\left(\frac{(\tau-\tau_*)^6}{\eps^2}\right).$$
Thus with Lemmas \ref{err_omega}, \ref{est_int} and \ref{cancFAR},
\begin{eqnarray*}
(\rm a)&&\int\limits_{t_-}^{t_l}\frac{\eps^kf(J,\ph,\tau)}{\omega^{2k-3}(\tau)}\,\dif t +\int\limits_{t_r}^{t_+}\frac{\eps^kf(J,\ph,\tau)}{\omega^{2k-3}(\tau)}\,\dif t\\
&=&\int\limits_{t_-}^{t_l}\frac{\eps^kf(J,\ph,\tau)}{\Omega^{2k-3}(\tau)}\,\dif t +\int\limits_{t_r}^{t_+}\frac{\eps^kf(J,\ph,\tau)}{\Omega^{2k-3}(\tau)}\,\dif t +\int\limits_{t_-}^{t_l}\frac{\eps^kO(1)}{\Omega^{2k-4}(\tau)}\,\dif t +\int\limits_{t_r}^{t_+}\frac{\eps^kO(1)}{\Omega^{2k-4}(\tau)}\,\dif t\\
&=&\int\limits_{t_-}^{t_l}\frac{\eps^kf(J,\ph,\tau)}{\Omega^{2k-3}(\tau)}\,\dif t +\int\limits_{t_r}^{t_+}\frac{\eps^kf(J,\ph,\tau)}{\Omega^{2k-3}(\tau)}\,\dif t +O(\eps^{\frac32})\\
&=&\int\limits_{t_-}^{t_l}\frac{\eps^kf(I_*,\ph_a,\tau_*)}{\Omega^{2k-3}(\tau)}\,\dif t +\int\limits_{t_-}^{t_l}\frac{\eps^kO(1)}{\Omega^{2k-4}(\tau)}\,\dif t +\int\limits_{t_-}^{t_l}\frac{\eps^{k-1}O(1)}{\Omega^{2k-6}(\tau)}\,\dif t +\int\limits_{t_-}^{t_l}\frac{\eps^{k-2}O(1)}{\Omega^{2k-9}(\tau)}\,\dif t\\
&&{}+\int\limits_{t_r}^{t_+}\frac{\eps^kf(I_*,\ph_a,\tau_*)}{\Omega^{2k-3}(\tau)}\,\dif t +\int\limits_{t_r}^{t_+}\frac{\eps^kO(1)}{\Omega^{2k-4}(\tau)}\,\dif t +\int\limits_{t_r}^{t_+}\frac{\eps^{k-1}O(1)}{\Omega^{2k-6}(\tau)}\,\dif t +\int\limits_{t_r}^{t_+}\frac{\eps^{k-2}O(1)}{\Omega^{2k-9}(\tau)}\,\dif t+O(\eps^{\frac32})\\
&=&O(\eps^{\frac32}).
\end{eqnarray*}
Applying the result of (a) and Lemma \ref{varbJ}, we obtain
\begin{eqnarray*}
(\rm b)&&\int\limits_{t_-}^{t_l}\frac{\eps^kf(J,\ph,\tau)}{\omega^{2k-2}(\tau)}\,\dif t +\int\limits_{t_r}^{t_+}\frac{\eps^kf(J,\ph,\tau)}{\omega^{2k-2}(\tau)}\,\dif t \\
&=&\int\limits_{t_-}^{t_l}\frac{\eps^kf(J,\ph,\tau)}{\Omega^{2k-2}(\tau)}\,\dif t +\int\limits_{t_-}^{t_l}\frac{\eps^k\tilde f(J,\ph,\tau)}{\Omega^{2k-3}(\tau)}\,\dif t  +\int\limits_{t_-}^{t_l}\frac{\eps^kO(1)}{\Omega^{2k-4}(\tau)}\,\dif t \\
&&{}+\int\limits_{t_r}^{t_+}\frac{\eps^kf(J,\ph,\tau)}{\Omega^{2k-2}(\tau)}\,\dif t +\int\limits_{t_r}^{t_+}\frac{\eps^k\tilde f(J,\ph,\tau)}{\Omega^{2k-3}(\tau)}\,\dif t +\int\limits_{t_r}^{t_+}\frac{\eps^kO(1)}{\Omega^{2k-4}(\tau)}\,\dif t\\
&=&\int\limits_{t_-}^{t_l}\frac{\eps^kf(J,\ph,\tau)}{\Omega^{2k-2}(\tau)}\,\dif t +\int\limits_{t_r}^{t_+}\frac{\eps^kf(J,\ph,\tau)}{\Omega^{2k-2}(\tau)}\,\dif t +O(\eps^{\frac32})\\
&=&\int\limits_{t_-}^{t_l}\frac{\eps^kf(I_*,\ph_a,\tau_*)}{\Omega^{2k-2}(\tau)}\,\dif t +\int\limits_{t_-}^{t_l}\frac{\eps^k\tilde{\tilde f}(I_*,\ph_a,\tau_*)}{\Omega^{2k-3}(\tau)}\,\dif t +\int\limits_{t_-}^{t_l}\frac{\eps^{k-1}\tilde{\tilde{\tilde f}}(I_*,\ph_a,\tau_*)}{\Omega^{2k-5}(\tau)}\,\dif t +\int\limits_{t_-}^{t_l}\frac{\eps^{k-2}O(1)}{\Omega^{2k-8}(\tau)}\,\dif t\\
&&{}+\int\limits_{t_r}^{t_+}\frac{\eps^kf(I_*,\ph_a,\tau_*)}{\Omega^{2k-2}(\tau)}\,\dif t +\int\limits_{t_r}^{t_+}\frac{\eps^k\tilde{\tilde f}(I_*,\ph_a,\tau_*)}{\Omega^{2k-3}(\tau)}\,\dif t +\int\limits_{t_r}^{t_+}\frac{\eps^{k-1}\tilde{\tilde{\tilde f}}(I_*,\ph_a,\tau_*)}{\Omega^{2k-5}(\tau)}\,\dif t\\
&&{}+\int\limits_{t_r}^{t_+}\frac{\eps^{k-2}O(1)}{\Omega^{2k-8}(\tau)}\,\dif t +O(\eps^{\frac32})\\
&=&\int\limits_{t_-}^{t_l}\frac{\eps^kf(I_*,\ph_a,\tau_*)}{\Omega^{2k-2}(\tau)}\,\dif t +\int\limits_{t_r}^{t_+}\frac{\eps^kf(I_*,\ph_a,\tau_*)}{\Omega^{2k-2}(\tau)}\,\dif t +O(\eps^{\frac32}).
\end{eqnarray*}
Here $\tilde f$, $\tilde{\tilde f}$ and $\tilde{\tilde{\tilde f}}$ are smooth functions.
\qed

\begin{lem}\label{J-Jpm}
For $t\in[t_-,t_l]$,
$$J-J_-=\int\limits_{t_-}^{t}\eps^5\frac{7{\omega'_*}^4\alpha^4(I_*,\ph_a,\tau_*)}{\Omega^8(\tau_1)}\,\dif t_1+O\left(\frac{\eps^4}{(\tau-\tau_*)^6}\right),$$
For $t\in[t_r,t_+]$, similarly we have
$$J-J_+=\int\limits_{t_+}^{t}\eps^5\frac{7{\omega'_*}^4\alpha^4(I_*,\ph_a,\tau_*)}{\Omega^8(\tau_1)}\,\dif t_1+O\left(\frac{\eps^4}{(\tau-\tau_*)^6}\right).$$
\end{lem}

\proof
Making use of of formula (\ref{dotJ}) and Lemma \ref{varbJ}, for $t\in[t_-,t_l]$, we get
\begin{eqnarray*}
&&J-J_-=\int\limits_{t_-}^{t}\dot J\,\dif t_1=\int\limits_{t_-}^{t}\eps^5\frac{7\big(\omega'(\tau_1)\big)^4\alpha^4(J,\ph,\tau_1)}{\omega^8(\tau_1)}\,\dif t_1+\int\limits_{t_-}^{t}\frac{\eps^5O(1)}{\omega^7}\,\dif t_1\\
&&\qquad\quad\ =\int\limits_{t_-}^{t}\eps^5\frac{7{\omega'_*}^4\alpha^4(I_*,\ph_a,\tau_*)}{\Omega^8(\tau_1)}\,\dif t_1+O\left(\frac{\eps^4}{(\tau-\tau_*)^6}\right).\\
\end{eqnarray*}
Similarly for $t\in[t_r,t_+]$, we obtain the second formula.

\subsubsection{Principal identity for formula (\ref{angle1})}

The relation between $\ph_+$ and $\ph_-$ is:
\begin{eqnarray*}
\ph_++\eps v_+&=&(\psi_+)\ =\ \psi_r+\int\limits_{t_r}^{t_+}\dot\psi\,\dif t\\
&=&\ph_r+\eps v_r+\int\limits_{t_r}^{t_+}\dot\psi\,\dif t\\
&=&\ph_l+\int\limits_{t_l}^{t_r}\omega(\tau)\,\dif t +\int\limits_{t_l}^{t_r}\eps\frac{\partial H_1(I,\ph,\tau)}{\partial I}\,\dif t +\eps v_r+\int\limits_{t_r}^{t_+}\dot\psi\,\dif t\\
&=&\psi_l-\eps v_l+\eps v_r+\int\limits_{t_l}^{t_r}\omega(\tau)\,\dif t +\int\limits_{t_l}^{t_r}\eps\frac{\partial H_1(I,\ph,\tau)}{\partial I}\,\dif t +\int\limits_{t_r}^{t_+}\dot\psi\,\dif t\\
&=&\ph_-+\eps v_-+\int\limits_{t_-}^{t_l}\dot\psi\,\dif t -\eps v_l+\eps v_r+\int\limits_{t_l}^{t_r}\omega(\tau)\,\dif t +\int\limits_{t_l}^{t_r}\eps\frac{\partial H_1(I,\ph,\tau)}{\partial I}\,\dif t +\int\limits_{t_r}^{t_+}\dot\psi\,\dif t\\
&=&\ph_-+\eps v_--\eps v_l+\eps v_r+\int\limits_{t_-}^{t_+}\omega(\tau)\,\dif t +\int\limits_{t_l}^{t_r}\eps\frac{\partial H_1(I,\ph,\tau)}{\partial I}\,\dif t\\
&&{}+\int\limits_{t_-}^{t_l}\frac{\partial}{\partial J}\left(\eps\RR_1+\eps^2\RR_2+\eps^3\RR_3+\eps^4\RR_4+\eps^5\dfrac{\partial S_4}{\partial\tau}\right)\,\dif t +\int\limits_{t_-}^{t_l}\sum\limits_{k=5}^{9}\eps^k\frac{\gamma^k(J,\ph,\tau)}{\omega^{2k-3}(\tau)}\,\dif t\\
&&{}+\int\limits_{t_r}^{t_+}\frac{\partial}{\partial J}\left(\eps\RR_1+\eps^2\RR_2+\eps^3\RR_3+\eps^4\RR_4+\eps^5\dfrac{\partial S_4}{\partial\tau}\right)\,\dif t +\int\limits_{t_r}^{t_+}\sum\limits_{k=5}^{9}\eps^k\frac{\gamma^k(J,\ph,\tau)}{\omega^{2k-3}(\tau)}\,\dif t+O(\eps^{\frac32}).
\end{eqnarray*}

Let
\begin{eqnarray*}
&&e_1^\ph=\eps\int\limits_{t_l}^{t_r}\frac{\partial H_1(I,\ph,\tau)}{\partial I}\,\dif t, \quad e_{j+1}^\ph=\int\limits_{t_-}^{t_l}\frac{\partial}{\partial J}\eps^j\RR_j\,\dif t+\int\limits_{t_r}^{t_+}\frac{\partial}{\partial J}\eps^j\RR_j\,\dif t, \quad j=1,2,3,4,\\
&&e_6^\ph=\int\limits_{t_-}^{t_l}\frac{\partial}{\partial J}\eps^5\dfrac{\partial S_4}{\partial\tau}\,\dif t+\int\limits_{t_r}^{t_+}\frac{\partial}{\partial J}\eps^5\dfrac{\partial S_4}{\partial\tau}\,\dif t,\quad e_7^\ph=\int\limits_{t_-}^{t_l}\sum\limits_{k=5}^{9}\eps^k\frac{\gamma^k(J,\ph,\tau)}{\omega^{2k-3}(\tau)}\,\dif t +\int\limits_{t_r}^{t_+}\sum\limits_{k=5}^{9}\eps^k\frac{\gamma^k(J,\ph,\tau)}{\omega^{2k-3}(\tau)}\,\dif t.
\end{eqnarray*}

The identity becomes
\begin{equation}\label{phiiden}
\ph_++\eps v_+=\ph_-+\eps v_-+\frac1\eps\int\limits_{\tau_-}^{\tau_+}\omega(\tau)\,\dif\tau-\eps v_l+\eps v_r+e_1^\ph+e_2^\ph+e_3^\ph+e_4^\ph+e_5^\ph+e_6^\ph+e_7^\ph.
\end{equation}
We will discuss the estimate term by term.

\subsubsection{Estimate of $-\eps v_l+\eps v_r$}

Here $v=v_1+v_2+v_3+v_4$, $v_k=\frac{\partial S_k}{\partial J}$, $k=1,2,3,4$. Making use of equations (\ref{dSdJ}) and Lemma \ref{replaceNEAR}, we get
\begin{eqnarray*}
-\eps v_l+\eps v_r&=&-\eps\frac{\partial S_l}{\partial J}+\eps\frac{\partial S_r}{\partial J}\\
&=&-\frac\eps\omega\frac{\partial\hat{\widetilde H}_1}{\partial J}\Bigg|_{\tau_l}^{\tau_r} -\eps^2\frac{\omega'}{\omega^3}\frac{\partial\hat{\hat{\widetilde H}}_1}{\partial J}\Bigg|_{\tau_l}^{\tau_r} -\eps^3\frac{3\omega'^2}{\omega^5}\frac{\partial\hat{\hat{\hat{\widetilde H}}}_1}{\partial J}\Bigg|_{\tau_l}^{\tau_r} -\eps^4\frac{15\omega'^3}{\omega^7}\frac{\partial\hat{\hat{\hat{\hat{\widetilde H}}}}_1}{\partial J}\Bigg|_{\tau_l}^{\tau_r} +\frac{\eps^2}{\omega^2}\frac{\partial\hat\beta^2}{\partial J}\Bigg|_{\tau_l}^{\tau_r} \\\\
&&{}+\frac{\eps^3}{\omega^4}\frac{\partial\hat\beta^3}{\partial J}\Bigg|_{\tau_l}^{\tau_r} +\frac{\eps^4}{\omega^6}\frac{\partial\hat\beta^4}{\partial J}\Bigg|_{\tau_l}^{\tau_r} \\\\
&=&-\frac\eps\Omega\frac{\partial\hat{\widetilde H}_1(I_*,\ph_a,\tau_*)}{\partial I}\Bigg|_{\tau_l}^{\tau_r} -\eps^2\frac{\Omega'}{\Omega^3}\frac{\partial\hat{\hat{\widetilde H}}_1(I_*,\ph_a,\tau_*)}{\partial I}\Bigg|_{\tau_l}^{\tau_r} -\eps^3\frac{3\Omega'^2}{\Omega^5}\frac{\partial\hat{\hat{\hat{\widetilde H}}}_1(I_*,\ph_a,\tau_*)}{\partial I}\Bigg|_{\tau_l}^{\tau_r}\\
&&{}-\eps^4\frac{15\Omega'^3}{\Omega^7}\frac{\partial\hat{\hat{\hat{\hat{\widetilde H}}}}_1(I_*,\ph_a,\tau_*)}{\partial I}\Bigg|_{\tau_l}^{\tau_r} +O(\eps^{\frac32}).
\end{eqnarray*}

\subsubsection{Estimate of $e_1^\ph$}

Using Lemmas \ref{cancNEAR} and \ref{varbI}, we have
\begin{eqnarray}\label{e1phi}
e_1^\ph&=&\int\limits_{t_l}^{t_r}\eps\frac{\partial H_1(I,\ph,\tau)}{\partial I}\,\dif t\nonumber\\
&=&\eps\int\limits_{t_l}^{t_r}\frac{\partial\bar H_1(I,\tau)}{\partial I}\,\dif t +\eps\int\limits_{t_l}^{t_r}\frac{\partial\widetilde H_1(I,\ph,\tau)}{\partial I}\,\dif t\nonumber\\
&=&\eps\int\limits_{t_l}^{t_*}\frac{\partial\bar H_1(I,\tau)}{\partial I}\,\dif t +\eps\int\limits_{t_*}^{t_r}\frac{\partial\bar H_1(I,\tau)}{\partial I}\,\dif t +\eps\int\limits_{t_l}^{t_r}\frac{\partial\widetilde H_1(I_*,\ph_a,\tau_*)}{\partial I}\,\dif t+O(\eps^{\frac32})\nonumber\\
&=&\eps\int\limits_{t_l}^{t_*}\frac{\partial\bar H_1(J_-,\tau)}{\partial I}\,\dif t +\eps\int\limits_{t_*}^{t_r}\frac{\partial\bar H_1(J_+,\tau)}{\partial I}\,\dif t +\eps\int\limits_{t_l}^{t_r}\frac{\partial\widetilde H_1(I_*,\ph_a,\tau_*)}{\partial I}\,\dif t\nonumber\\
&&{}+\eps\int\limits_{t_l}^{t_*}\frac{\partial^2\bar H_1(J_-,\tau)}{\partial I^2}(I-J_-)\,\dif t+\eps\int\limits_{t_*}^{t_r}\frac{\partial^2\bar H_1(J_+,\tau)}{\partial I^2}(I-J_+)\,\dif t+O(\eps^{\frac32})\nonumber\\
&=&\eps\int\limits_{t_l}^{t_*}\frac{\partial\bar H_1(J_-,\tau)}{\partial I}\,\dif t +\eps\int\limits_{t_*}^{t_r}\frac{\partial\bar H_1(J_+,\tau)}{\partial I}\,\dif t +\eps\int\limits_{t_l}^{t_r}\frac{\partial\widetilde H_1(I_*,\ph_a,\tau_*)}{\partial I}\,\dif t\nonumber\\
&&{}+\eps\int\limits_{t_l}^{t_*}\frac{\partial^2\bar H_1(I_*,\tau_*)}{\partial I^2}(I-J_-)\,\dif t+\eps\int\limits_{t_*}^{t_r}\frac{\partial^2\bar H_1(I_*,\tau_*)}{\partial I^2}(I-J_+)\,\dif t+O(\eps^{\frac32}).
\end{eqnarray}

\begin{lem}\label{cancIJ+-}
$$\eps\int\limits_{t_l}^{t_*}(I-J_-)\,\dif t+\eps\int\limits_{t_*}^{t_r}(I-J_+)\,\dif t=O(\eps^{\frac32}).$$
\end{lem}

\proof
Making use of Lemmas \ref{varbI}, \ref{varbJ} and \ref{J-Jpm}, we get
\begin{eqnarray*}
&&\eps\int\limits_{t_l}^{t_*}(I-J_-)\,\dif t+\eps\int\limits_{t_*}^{t_r}(I-J_+)\,\dif t\\
&=&\eps\int\limits_{t_l}^{t_*}(I-I_l+I_l-J_l+J_l-J_-)\,\dif t +\eps\int\limits_{t_*}^{t_r}(I-I_r+I_r-J_r+J_r-J_+)\,\dif t\\
&=&\eps\int\limits_{t_l}^{t_*}\left(-\eps\int\limits_{t_l}^{t}\frac{\partial H_1(I,\ph,\tau_1)}{\partial\ph}\,\dif t_1+\eps\frac{\partial S_l}{\partial\ph}+\int\limits_{t_-}^{t_l}\dot J\,\dif t\right)\,\dif t \\
&&{}+\eps\int\limits_{t_*}^{t_r}\left(\eps\int\limits_{t}^{t_r}\frac{\partial H_1(I,\ph,\tau_1)}{\partial\ph}\,\dif t_1+\eps\frac{\partial S_r}{\partial\ph}+\int\limits_{t_+}^{t_r}\dot J\,\dif t\right)\,\dif t \\
&=&\eps\int\limits_{t_l}^{t_*}\left[-\eps\int\limits_{t_l}^{t}\frac{\partial H_1(I_*,\ph_a,\tau_*)}{\partial\ph}\,\dif t_1+\eps\frac{\alpha^1(J_l,\ph_l,\tau_l)}{\omega(\tau_l)}+\eps^2\frac{\omega'(\tau_l)\alpha^2(J_l,\ph_l,\tau_l)}{\omega^3(\tau_l)}\right.\\
&&\qquad{}+\left.\eps^3\frac{\big(\omega'(\tau_l)\big)^2\alpha^3(J_l,\ph_l,\tau_l)}{\omega^5(\tau_l)}+\eps^4\frac{\big(\omega'(\tau_l)\big)^3\alpha^4(J_l,\ph_l,\tau_l)}{\omega^7(\tau_l)}+\int\limits_{t_-}^{t_l}\dot J\,\dif t\right]\,\dif t \\
&&{}+\eps\int\limits_{t_*}^{t_r}\left[\eps\int\limits_{t}^{t_r}\frac{\partial H_1(I_*,\ph_a,\tau_*)}{\partial\ph}\,\dif t_1+\eps\frac{\alpha^1(J_r,\ph_r,\tau_r)}{\omega(\tau_r)}+\eps^2\frac{\omega'(\tau_r)\alpha^2(J_r,\ph_r,\tau_r)}{\omega^3(\tau_r)}\right.\\
&&\qquad{}+\eps^3\frac{\big(\omega'(\tau_r)\big)^2\alpha^3(J_r,\ph_r,\tau_r)}{\omega^5(\tau_r)}+\eps^4\frac{\big(\omega'(\tau_r)\big)^3\alpha^4(J_r,\ph_r,\tau_r)}{\omega^7(\tau_r)}+\left.\int\limits_{t_+}^{t_r}\dot J\,\dif t\right]\,\dif t +O(\eps^{\frac32})\\
&=&\eps\int\limits_{t_l}^{t_*}\Bigg[\eps\frac{\alpha^1(I_*,\ph_{a_l},\tau_*)}{\Omega(\tau_l)}+\eps^2\frac{\Omega'_*\alpha^2(I_*,\ph_{a_l},\tau_*)}{\Omega^3(\tau_l)}+\eps^3\frac{{\Omega'_*}^2\alpha^3(I_*,\ph_{a_l},\tau_*)}{\Omega^5(\tau_l)}+\eps^4\frac{{\Omega'_*}^3\alpha^4(I_*,\ph_{a_l},\tau_*)}{\Omega^7(\tau_l)}\\
&&\qquad{}+\int\limits_{t_-}^{t_l}\eps^5\frac{7{\omega'_*}^4\alpha^4(I_*,\ph_a,\tau_*)}{\Omega^8(\tau_1)}\,\dif t_1 +\eps\int\limits_{t_*}^{t_r}\Bigg[\eps\frac{\alpha^1(I_*,\ph_{a_r},\tau_*)}{\Omega(\tau_r)}+\eps^2\frac{\Omega'_*\alpha^2(I_*,\ph_{a_r},\tau_*)}{\Omega^3(\tau_r)}\\
&&\qquad{}+\eps^3\frac{{\Omega'_*}^2\alpha^3(I_*,\ph_{a_r},\tau_*)}{\Omega^5(\tau_r)}+\eps^4\frac{{\Omega'_*}^3\alpha^4(I_*,\ph_{a_r},\tau_*)}{\Omega^7(\tau_r)}+\int\limits_{t_+}^{t_r}\eps^5\frac{7{\omega'_*}^4\alpha^4(I_*,\ph_a,\tau_*)}{\Omega^8(\tau_1)}\,\dif t_1\,\Bigg]\,\dif t+O(\eps^{\frac32})\\
&=&O(\eps^{\frac32}).
\end{eqnarray*}
\qed

We also have (cf. Section \ref{e4}):
\begin{eqnarray*}
&&\eps\int\limits_{t_l}^{t_r}\frac{\partial\widetilde H_1(I_*,\ph_a,\tau_*)}{\partial I}\,\dif t\\
&=&\eps\int\limits_{-\infty}^{+\infty}\frac{\partial\widetilde H_1(I_*,\ph_a,\tau_*)}{\partial I}\,\dif t-\eps\int\limits_{t_-}^{t_l}\frac{\partial\widetilde H_1(I_*,\ph_a,\tau_*)}{\partial I}\,\dif t -\eps\int\limits_{t_r}^{t_+}\frac{\partial\widetilde H_1(I_*,\ph_a,\tau_*)}{\partial I}\,\dif t\\
&&{}-\eps\int\limits_{-\infty}^{t_-}\frac{\partial\widetilde H_1(I_*,\ph_a,\tau_*)}{\partial I}\,\dif t-\eps\int\limits_{t_+}^{+\infty}\frac{\partial\widetilde H_1(I_*,\ph_a,\tau_*)}{\partial I}\,\dif t\\
&=&\eps\int\limits_{-\infty}^{+\infty}\frac{\partial\widetilde H_1(I_*,\ph_a,\tau_*)}{\partial I}\,\dif t-\eps\int\limits_{t_-}^{t_l}\frac{\partial\widetilde H_1(I_*,\ph_a,\tau_*)}{\partial I}\,\dif t -\eps\int\limits_{t_r}^{t_+}\frac{\partial\widetilde H_1(I_*,\ph_a,\tau_*)}{\partial I}\,\dif t\\
&&{}-\frac{\eps}{\Omega(\tau_-)}\frac{\partial \hat{\widetilde H}_1(I_*,\ph_{a_-},\tau_*)}{\partial I}+\frac{\eps}{\Omega(\tau_+)}\frac{\partial \hat{\widetilde H}_1(I_*,\ph_{a_+},\tau_*)}{\partial I}+O(\eps^2)\\
&=&\sqrt\eps\int\limits_{-\infty}^{+\infty}\frac{\partial\widetilde H_1(I_*,\ph_*+\frac{\omega'_*}{2}\theta^2,\tau_*)}{\partial I}\,\dif\theta -\eps\int\limits_{t_-}^{t_l}\frac{\partial\widetilde H_1(I_*,\ph_a,\tau_*)}{\partial I}\,\dif t -\eps\int\limits_{t_r}^{t_+}\frac{\partial\widetilde H_1(I_*,\ph_a,\tau_*)}{\partial I}\,\dif t\\
&&{}+\frac{\eps}{\Omega(\tau)}\frac{\partial \hat{\widetilde H}_1(I_*,\ph_a,\tau_*)}{\partial I}\Bigg|_{\tau_-}^{\tau_+}+O(\eps^2).\\
\end{eqnarray*}

Using this result, formula (\ref{e1phi}) and Lemma \ref{cancIJ+-}, we obtain
\begin{eqnarray*}
e_1^\ph&=&\eps\int\limits_{t_l}^{t_*}\frac{\partial\bar H_1(J_-,\tau)}{\partial I}\,\dif t +\eps\int\limits_{t_*}^{t_r}\frac{\partial\bar H_1(J_+,\tau)}{\partial I}\,\dif t +\sqrt\eps\int\limits_{-\infty}^{+\infty}\frac{\partial\widetilde H_1(I_*,\ph_*+\frac{\omega'_*}{2}\theta^2,\tau_*)}{\partial I}\,\dif\theta\\ &&-\eps\int\limits_{t_-}^{t_l}\frac{\partial\widetilde H_1(I_*,\ph_a,\tau_*)}{\partial I}\,\dif t -\eps\int\limits_{t_r}^{t_+}\frac{\partial\widetilde H_1(I_*,\ph_a,\tau_*)}{\partial I}\,\dif t+\frac{\eps}{\Omega}\frac{\partial \hat{\widetilde H}_1(I_*,\ph_a,\tau_*)}{\partial I}\Bigg|_{t_-}^{t_+}+O(\eps^{\frac32}).
\end{eqnarray*}

\subsubsection{Estimate of $e_2^\ph$}

We know that $\RR_1=\big\langle H_1\big\rangle^\ph=\bar H_1$. Thus with Lemma \ref{J-Jpm}, we have
\begin{eqnarray*}
e_2^\ph&=&\int\limits_{t_-}^{t_l}\frac{\partial}{\partial J}\,\eps\RR_1\,\dif t +\int\limits_{t_r}^{t_+}\frac{\partial}{\partial J}\,\eps\RR_1\,\dif t\\
&=&\eps\int\limits_{t_-}^{t_l}\frac{\partial\bar H_1(J,\tau)}{\partial J}\,\dif t +\eps\int\limits_{t_r}^{t_+}\frac{\partial\bar H_1(J,\tau)}{\partial J}\,\dif t\\
&=&\eps\int\limits_{t_-}^{t_l}\frac{\partial\bar H_1(J_-,\tau)}{\partial I}\,\dif t +\eps\int\limits_{t_-}^{t_l}\frac{\partial^2\bar H_1(J_-,\tau)}{\partial I^2}(J-J_-)\,\dif t +\eps\int\limits_{t_-}^{t_l}O(J-J_-)^2\,\dif t \\
&&{}+\eps\int\limits_{t_r}^{t_+}\frac{\partial\bar H_1(J_+,\tau)}{\partial I}\,\dif t +\eps\int\limits_{t_r}^{t_+}\frac{\partial^2\bar H_1(J_+,\tau)}{\partial I^2}(J-J_+)\,\dif t +\eps\int\limits_{t_r}^{t_+}O(J-J_+)^2\,\dif t\\
&=&\eps\int\limits_{t_-}^{t_l}\frac{\partial\bar H_1(J_-,\tau)}{\partial I}\,\dif t +\eps\int\limits_{t_-}^{t_l}\frac{\partial^2\bar H_1(I_*,\tau_*)}{\partial I^2}\int\limits_{t_-}^{t}\eps^5\frac{7{\omega'_*}^4\alpha^4(I_*,\ph_a,\tau_*)}{\Omega^8(\tau_1)}\,\dif t_1\,\dif t\\
&&{}+\eps\int\limits_{t_r}^{t_+}\frac{\partial\bar H_1(J_+,\tau)}{\partial I}\,\dif t +\eps\int\limits_{t_r}^{t_+}\frac{\partial^2\bar H_1(I_*,\tau_*)}{\partial I^2}\int\limits_{t_+}^{t}\eps^5\frac{7{\omega'_*}^4\alpha^4(I_*,\ph_a,\tau_*)}{\Omega^8(\tau_1)}\,\dif t_1\,\dif t +O(\eps^{\frac32})\\
&=&\eps\int\limits_{t_-}^{t_l}\frac{\partial\bar H_1(J_-,\tau)}{\partial I}\,\dif t +\eps\int\limits_{t_r}^{t_+}\frac{\partial\bar H_1(J_+,\tau)}{\partial I}\,\dif t +O(\eps^{\frac32}).
\end{eqnarray*}

\subsubsection{Estimate of $e_3^\ph$}

From Lemma \ref{varbJ}, we get the error between $J$ and $I_*$:
\begin{eqnarray*}
J-I_*&=&-\int\limits_{\tau_*}^{\tau}\frac{\partial H_1(I_*,\ph_a,\tau_*)}{\partial\ph}\,\dif\tau_1+\eps\frac{f^1(I_*,\ph_a,\tau_*)}{\tau-\tau_*}+\eps^2\frac{f^2(I_*,\ph_a,\tau_*)}{(\tau-\tau_*)^3}+\eps^3\frac{f^3(I_*,\ph_a,\tau_*)}{(\tau-\tau_*)^5}\\
&&{}+\eps^4\frac{f^4(I_*,\ph_a,\tau_*)}{(\tau-\tau_*)^7}+O(\eps)+O(\tau-\tau_*)^2+O\left(\frac{(\tau-\tau_*)^4}{\eps}\right)+O\left(\frac{(\tau-\tau_*)^6}{\eps^2}\right).
\end{eqnarray*}

As
\begin{eqnarray*}
\RR_2&=&\left\langle\frac{\partial H_1}{\partial I}\frac{\partial S_1}{\partial\ph}\right\rangle^\ph=\left\langle\frac{\partial H_1}{\partial I}\frac{\alpha^1}{\omega}\right\rangle^\ph=\left\langle\frac{\partial\widetilde H_1}{\partial I}\cdot\left(\frac{-\widetilde H_1}{\omega}\right)\right\rangle^\ph\\
&=&-\frac1\omega\left\langle\frac12\frac{\partial}{\partial I}\widetilde H_1^2\right\rangle^\ph=-\frac1{2\omega}\frac{\partial}{\partial I}\big\langle\widetilde H_1^2\big\rangle^\ph,
\end{eqnarray*}

we have
\begin{eqnarray*}
e_3^\ph&=&\int\limits_{t_-}^{t_l}\frac{\partial}{\partial J}\,\eps^2\RR_2\,\dif t +\int\limits_{t_r}^{t_+}\frac{\partial}{\partial J}\,\eps^2\RR_2\,\dif t\\
&=&\eps^2\int\limits_{t_-}^{t_l}\left(-\frac1{2\omega}\right) \frac{\partial^2\big\langle\widetilde H_1^2\big\rangle^\ph(J,\tau)}{\partial I^2}\,\dif t+\eps^2\int\limits_{t_r}^{t_+}\left(-\frac1{2\omega}\right) \frac{\partial^2\big\langle\widetilde H_1^2\big\rangle^\ph(J,\tau)}{\partial I^2}\,\dif t\\
&=&\eps^2\int\limits_{t_-}^{t_l}\left(-\frac1{2\omega}\right) \frac{\partial^2\big\langle\widetilde H_1^2\big\rangle^\ph(J_-,\tau)}{\partial I^2}\,\dif t+\eps^2\int\limits_{t_-}^{t_l}\left(-\frac1{2\omega}\right) \frac{\partial^3\big\langle\widetilde H_1^2\big\rangle^\ph(J_-,\tau)}{\partial I^3}(J-J_-)\,\dif t\\
&&{}+\eps^2\int\limits_{t_r}^{t_+}\left(-\frac1{2\omega}\right) \frac{\partial^2\big\langle\widetilde H_1^2\big\rangle^\ph(J_+,\tau)}{\partial I^2}\,\dif t+\eps^2\int\limits_{t_r}^{t_+}\left(-\frac1{2\omega}\right) \frac{\partial^3\big\langle\widetilde H_1^2\big\rangle^\ph(J_+,\tau)}{\partial I^3}(J-J_+)\,\dif t+O(\eps^2)\\
&=&\eps^2\int\limits_{t_-}^{t_l}\left(-\frac1{2\omega}\right) \frac{\partial^2\big\langle\widetilde H_1^2\big\rangle^\ph(J_-,\tau)}{\partial I^2}\,\dif t+\eps^2\int\limits_{t_-}^{t_l}\left(-\frac1{2\omega}\right) \frac{\partial^3\big\langle\widetilde H_1^2\big\rangle^\ph(I_*,\tau_*)}{\partial I^3}\int\limits_{t_-}^{t}\eps^5\frac{7{\omega'_*}^4\alpha^4(I_*,\ph_a,\tau_*)}{\Omega^8(\tau_1)}\,\dif t_1\,\dif t\\
&&{}+\eps^2\int\limits_{t_r}^{t_+}\left(-\frac1{2\omega}\right) \frac{\partial^2\big\langle\widetilde H_1^2\big\rangle^\ph(J_+,\tau)}{\partial I^2}\,\dif t+\eps^2\int\limits_{t_r}^{t_+}\left(-\frac1{2\omega}\right) \frac{\partial^3\big\langle\widetilde H_1^2\big\rangle^\ph(I_*,\tau_*)}{\partial I^3}\int\limits_{t_+}^{t}\eps^5\frac{7{\omega'_*}^4\alpha^4(I_*,\ph_a,\tau_*)}{\Omega^8(\tau_1)}\,\dif t_1\,\dif t\\
&&{}+O(\eps^2)\\
&=&\eps^2\int\limits_{t_-}^{t_l}\left(-\frac1{2\omega}\right) \frac{\partial^2\big\langle\widetilde H_1^2\big\rangle^\ph(J_-,\tau)}{\partial I^2}\,\dif t+\eps^2\int\limits_{t_r}^{t_+}\left(-\frac1{2\omega}\right) \frac{\partial^2\big\langle\widetilde H_1^2\big\rangle^\ph(J_+,\tau)}{\partial I^2}\,\dif t+O(\eps^\frac32)\\
&=&\eps^2\int\limits_{t_-}^{t_l}\left(-\frac1{2\omega}\right) \frac{\partial^2\big\langle\widetilde H_1^2\big\rangle^\ph(J_-,\tau)}{\partial I^2}\,\dif t+\eps^2\int\limits_{t_r}^{t_+}\left(-\frac1{2\omega}\right) \frac{\partial^2\big\langle\widetilde H_1^2\big\rangle^\ph(J_-,\tau)}{\partial I^2}\,\dif t\\
&&{}+\eps^2\int\limits_{t_r}^{t_+}\left(-\frac1{2\omega}\right) \frac{\partial^3\big\langle\widetilde H_1^2\big\rangle^\ph(J_-,\tau)}{\partial I^3}(J_+-J_-)\,\dif t +O(\eps^\frac32)\\
\end{eqnarray*}
\begin{eqnarray*}
&=&{\rm p.v.}\,\eps^2\int\limits_{t_-}^{t_+}\frac{\partial{\cal R}_2(J_-,\tau)}{\partial I}\,\dif t -\lim_{\delta\to0+}\left[\eps^2\int\limits_{t_l}^{t_*-\delta} \left(-\frac1{2\omega}\right) \frac{\partial^2\big\langle\widetilde H_1^2\big\rangle^\ph}{\partial I^2}\,\dif t +\eps^2\int\limits_{t_*+\delta}^{t_r}\left(-\frac1{2\omega}\right) \frac{\partial^2\big\langle\widetilde H_1^2\big\rangle^\ph}{\partial I^2}\,\dif t \right]\\
&&{}+\eps^2\int\limits_{t_r}^{t_+}\left(-\frac1{2\omega}\right) \frac{\partial^3\big\langle\widetilde H_1^2\big\rangle^\ph(J_-,\tau_*)}{\partial I^3} \left(\,\int\limits_{-\infty}^{+\infty}\dot I_a\,\dif t+O(\eps^{\frac32})\right)\,\dif t +O(\eps^{\frac32})\\
&=&{\rm p.v.}\,\eps^2\int\limits_{t_-}^{t_+}\frac{\partial{\RR}_2(J_-,\tau)}{\partial I}\,\dif t+\eps^3\frac{\partial^3\big\langle\widetilde H_1^2\big\rangle^\ph(J_-,\tau_*)}{\partial I^3}\int\limits_{-\infty}^{+\infty}\frac{\partial H_1(I_*,\ph_a,\tau_*)}{\partial\ph}\,\dif t\int\limits_{t_r}^{t_+}\frac1{2\omega}\,\dif t +O(\eps^{\frac32})\\
&=&{\rm p.v.}\,\eps^2\int\limits_{t_-}^{t_+}\frac{\partial{\RR}_2(I_*,\tau)}{\partial I}\,\dif t-\frac{\eps^2\ln\eps}{4\omega'_*}\frac{\partial^3\big\langle\widetilde H_1^2\big\rangle^\ph(I_*,\tau_*)}{\partial I^3} \int\limits_{-\infty}^{+\infty}\frac{\partial H_1(I_*,\ph_a,\tau_*)}{\partial\ph}\,\dif t+O(\eps^{\frac32}).
\end{eqnarray*}

\subsubsection{Estimate of $e_4^\ph$}

\begin{eqnarray*}
\RR_3&=&\left\langle\frac{\partial H_1}{\partial I}\frac{\partial S_2}{\partial\ph} +\frac12\frac{\partial^2H_1}{\partial I^2}\left(\frac{\partial S_1}{\partial\ph}\right)^2 \right\rangle^\ph\\
&=&\left\langle\frac{\partial H_1}{\partial I}\frac{-\omega'\hat{\widetilde H}_1}{\omega^3} +\frac12\frac{O(1)}{\omega^2} \right\rangle^\ph\\
&=&-\frac{\omega'}{\omega^3}\left\langle\hat{\widetilde H}_1\frac{\partial H_1}{\partial I}\right\rangle^\ph+\frac{O(1)}{\omega^2}.
\end{eqnarray*}

So
\begin{eqnarray*}
e_4^\ph&=&\int\limits_{t_-}^{t_l}\frac{\partial}{\partial J}\,\eps^3\RR_3\,\dif t +\int\limits_{t_r}^{t_+}\frac{\partial}{\partial J}\,\eps^3\RR_3\,\dif t\\
&=&-\eps^3\int\limits_{t_-}^{t_l}\frac{\omega'}{\omega^3}\frac{\partial}{\partial J} \left\langle\hat{\widetilde H}_1\frac{\partial H_1}{\partial I}\right\rangle^\ph(J,\tau)\,\dif t-\eps^3\int\limits_{t_r}^{t_+}\frac{\omega'}{\omega^3}\frac{\partial}{\partial J} \left\langle\hat{\widetilde H}_1\frac{\partial H_1}{\partial I}\right\rangle^\ph(J,\tau)\,\dif t+O(\eps^{\frac32})\\
&=&-\eps^3\int\limits_{t_-}^{t_l}\frac{\Omega'}{\Omega^3}\frac{\partial}{\partial J} \left\langle\hat{\widetilde H}_1\frac{\partial H_1}{\partial I}\right\rangle^\ph(I_*,\tau_*)\,\dif t-\eps^3\int\limits_{t_r}^{t_+}\frac{\Omega'}{\Omega^3}\frac{\partial}{\partial J} \left\langle\hat{\widetilde H}_1\frac{\partial H_1}{\partial I}\right\rangle^\ph(I_*,\tau_*)\,\dif t+O(\eps^{\frac32})\\
&=&O(\eps^{\frac32}).
\end{eqnarray*}

\subsubsection{Estimate of $e_5^\ph$}

\begin{eqnarray*}
\RR_4&=&\left\langle\frac{\partial H_1}{\partial I}\frac{\partial S_3}{\partial\ph} +\frac{\partial^2H_1}{\partial I^2}\frac{\partial S_1}{\partial\ph}\frac{\partial S_2}{\partial\ph}+\frac16\frac{\partial^3H_1}{\partial I^3}\left(\frac{\partial S_1}{\partial\ph}\right)^{\!3}\,\right\rangle^\ph\\
&=&\left\langle\frac{\partial H_1}{\partial I}\frac{-3\omega'^2\hat{\hat{\widetilde H}}_1}{\omega^5} +\frac{O(1)}{\omega^4} +\frac16\frac{O(1)}{\omega^3} \right\rangle^\ph\\
&=&-\frac{3\omega'^2}{\omega^5}\left\langle\hat{\hat{\widetilde H}}_1\frac{\partial H_1}{\partial I}\right\rangle^\ph+\frac{O(1)}{\omega^4}.
\end{eqnarray*}

So
\begin{eqnarray*}
e_5^\ph&=&\int\limits_{t_-}^{t_l}\frac{\partial}{\partial J}\,\eps^4\RR_4\,\dif t +\int\limits_{t_r}^{t_+}\frac{\partial}{\partial J}\,\eps^4\RR_4\,\dif t\\
&=&-\eps^4\int\limits_{t_-}^{t_l}\frac{3\omega'^2}{\omega^5}\frac{\partial}{\partial J} \left\langle\hat{\hat{\widetilde H}}_1\frac{\partial H_1}{\partial I}\right\rangle^\ph(J,\tau)\,\dif t-\eps^4\int\limits_{t_r}^{t_+}\frac{3\omega'^2}{\omega^5}\frac{\partial}{\partial J} \left\langle\hat{\hat{\widetilde H}}_1\frac{\partial H_1}{\partial I}\right\rangle^\ph(J,\tau)\,\dif t+O(\eps^{\frac32})\\
&=&-\eps^4\int\limits_{t_-}^{t_l}\frac{3\Omega'^2}{\Omega^5}\frac{\partial}{\partial J} \left\langle\hat{\hat{\widetilde H}}_1\frac{\partial H_1}{\partial I}\right\rangle^\ph(I_*,\tau_*)\,\dif t-\eps^4\int\limits_{t_r}^{t_+}\frac{3\Omega'^2}{\Omega^5}\frac{\partial}{\partial J} \left\langle\hat{\hat{\widetilde H}}_1\frac{\partial H_1}{\partial I}\right\rangle^\ph(I_*,\tau_*)\,\dif t+O(\eps^{\frac32})\\
&=&O(\eps^{\frac32}).
\end{eqnarray*}

\subsubsection{Estimate of $e_6^\ph$}

There exists a function $U(I_*,\ph_a,\tau_*)$, such that $\dfrac{\partial U}{\partial\ph}=\dfrac{\partial \widetilde H_1(I_*,\ph_a,\tau_*)}{\partial I}$. Here\\ $U=\dfrac{\partial \hat{\widetilde H}_1(I_*,\ph_a,\tau_*)}{\partial I}$ with $\left\langle U \right\rangle^{\ph_a}=0$. Thus

\begin{eqnarray*}
\eps\int\frac{\partial \widetilde H_1(I_*,\ph_a,\tau_*)}{\partial I}\,\dif t &=&\eps\int\frac{\partial U}{\partial\ph}\,\dif t =\eps\int\frac{\dif U}{\dif t}\cdot\frac1\Omega\,\dif t =\eps\int\frac1\Omega\,\dif U \\
&=&\frac{\eps U(I_*,\ph_a,\tau_*)}{\Omega(\tau)}+\eps^2\int U\frac{\Omega'}{\Omega^2}\,\dif t.
\end{eqnarray*}

There exists a function $V(I_*,\ph_a,\tau_*)$, such that $\dfrac{\partial V}{\partial\ph}=U(I_*,\ph_a,\tau_*)$. Here\\ $V=\hat U(I_*,\ph_a,\tau_*)=\dfrac{\partial \hat{\hat{\widetilde H}}_1(I_*,\ph_a,\tau_*)}{\partial I}$ with $\left\langle V \right\rangle^{\ph_a}=0$. Thus

\begin{eqnarray*}
\eps^2\int U\frac{\Omega'}{\Omega^2}\,\dif t&=&\eps^2\int\frac{\partial V}{\partial\ph}\frac{\Omega'}{\Omega^2}\,\dif t=\eps^2\int\frac{\dif V}{\dif t}\frac{\Omega'}{\Omega^3}\,\dif t=\eps^2\int\frac{\Omega'}{\Omega^3}\,\dif V\\
&=&\eps^2V(I_*,\ph_a,\tau_*)\frac{\Omega'(\tau)}{\Omega^3(\tau)}+\eps^3\int V\frac{3\Omega'^2}{\Omega^4}\,\dif t.
\end{eqnarray*}

There exists a function $W(I_*,\ph_a,\tau_*)$, such that $\dfrac{\partial W}{\partial\ph}=V(I_*,\ph_a,\tau_*)$. Here\\ $W=\hat V(I_*,\ph_a,\tau_*)=\dfrac{\partial \hat{\hat{\hat{\widetilde H}}}_1(I_*,\ph_a,\tau_*)}{\partial I}$ with $\left\langle W \right\rangle^{\ph_a}=0$. Thus

\begin{eqnarray*}
\eps^3\int V\frac{3\Omega'^2}{\Omega^4}\,\dif t&=&\eps^3\int\frac{\partial W}{\partial\ph}\frac{3\Omega'^2}{\Omega^4}\,\dif t=\eps^3\int\frac{\dif W}{\dif t}\frac{3\Omega'^2}{\Omega^5}\,\dif t=\eps^3\int\frac{3\Omega'^2}{\Omega^5}\,\dif W\\
&=&\eps^3W(I_*,\ph_a,\tau_*)\frac{3\Omega'^2(\tau)}{\Omega^5(\tau)}+\eps^4\int W\frac{15\Omega'^3}{\Omega^6}\,\dif t.
\end{eqnarray*}

There exists a function $T(I_*,\ph_a,\tau_*)$, such that $\dfrac{\partial T}{\partial\ph}=W(I_*,\ph_a,\tau_*)$. Here\\ $T=\hat W(I_*,\ph_a,\tau_*)=\dfrac{\partial \hat{\hat{\hat{\hat{\widetilde H}}}}_1(I_*,\ph_a,\tau_*)}{\partial I}$ with $\left\langle T \right\rangle^{\ph_a}=0$. Thus

\begin{eqnarray*}
\eps^4\int W\frac{15\Omega'^3}{\Omega^6}\,\dif t&=&\eps^4\int\frac{\partial T}{\partial\ph}\frac{15\Omega'^3}{\Omega^6}\,\dif t=\eps^4\int\frac{\dif T}{\dif t}\frac{15\Omega'^3}{\Omega^7}\,\dif t=\eps^4\int\frac{15\Omega'^3}{\Omega^7}\,\dif T\\
&=&\eps^4T(I_*,\ph_a,\tau_*)\frac{15\Omega'^3(\tau)}{\Omega^7(\tau)}+\eps^5\int T\frac{105\Omega'^4}{\Omega^8}\,\dif t.
\end{eqnarray*}

Thus we have
\begin{eqnarray*}
\eps\int\frac{\partial \widetilde H_1(I_*,\ph_a,\tau_*)}{\partial I}\,\dif t &=&\frac{\eps U}{\Omega}+\eps^2V\frac{\Omega'}{\Omega^3}+\eps^3W\frac{3\Omega'^2}{\Omega^5}+\eps^4T\frac{15\Omega'^3}{\Omega^7}+\eps^5\int T\frac{105\Omega'^4}{\Omega^8}\,\dif t\\\\
&=&\frac{\eps}{\Omega}\dfrac{\partial \hat{\widetilde H}_1}{\partial I} +\eps^2\frac{\Omega'}{\Omega^3}\dfrac{\partial \hat{\hat{\widetilde H}}_1}{\partial I} +\eps^3\frac{3\Omega'^2}{\Omega^5}\dfrac{\partial \hat{\hat{\hat{\widetilde H}}}_1}{\partial I}+\eps^4\frac{15\Omega'^3}{\Omega^7}\dfrac{\partial \hat{\hat{\hat{\hat{\widetilde H}}}}_1}{\partial I}+\eps^5\int \frac{105\Omega'^4}{\Omega^8}\dfrac{\partial \hat{\hat{\hat{\hat{\widetilde H}}}}_1}{\partial I}\,\dif t.
\end{eqnarray*}

Therefore, applying equation (\ref{S}) and Lemma \ref{replaceFAR}, we obtain
\begin{eqnarray*}
e_6^\ph&=&\int\limits_{t_-}^{t_l}\frac{\partial}{\partial J}\,\eps^5\frac{\partial S_4}{\partial\tau}\,\dif t +\int\limits_{t_r}^{t_+}\frac{\partial}{\partial J}\,\eps^5\frac{\partial S_4}{\partial\tau}\,\dif t\\
&=&\int\limits_{t_-}^{t_l}\eps^5\frac{105\omega'^4}{\omega^8}\frac{\partial \hat{\hat{\hat{\hat{\widetilde H}}}}_1(J,\ph,\tau)}{\partial J}\,\dif t +\int\limits_{t_-}^{t_l}\frac{\eps^5\widetilde\beta(J,\ph,\tau)}{\omega^7}\,\dif t\\
&&{}+\int\limits_{t_r}^{t_+}\eps^5\frac{105\omega'^4}{\omega^8}\frac{\partial \hat{\hat{\hat{\hat{\widetilde H}}}}_1(J,\ph,\tau)}{\partial J}\,\dif t +\int\limits_{t_r}^{t_+}\frac{\eps^5\widetilde\beta(J,\ph,\tau)}{\omega^7}\,\dif t\\
&=&\eps^5\int\limits_{t_-}^{t_l}\frac{105\Omega'^4}{\Omega^8}\frac{\partial \hat{\hat{\hat{\hat{\widetilde H}}}}_1(I_*,\ph_a,\tau_*)}{\partial I}\,\dif t +\eps^5\int\limits_{t_r}^{t_+}\frac{105\Omega'^4}{\Omega^8}\frac{\partial \hat{\hat{\hat{\hat{\widetilde H}}}}_1(I_*,\ph_a,\tau_*)}{\partial I}\,\dif t +O(\eps^{\frac32})\\
&=&\eps\int\limits_{t_-}^{t_l}\frac{\partial \widetilde H_1(I_*,\ph_a,\tau_*)}{\partial I}\,\dif t -\frac{\eps}{\Omega(\tau)}\dfrac{\partial \hat{\widetilde H}_1(I_*,\ph_a,\tau_*)}{\partial I}\Bigg|_{t_-}^{t_l} -\eps^2\frac{\Omega'(\tau)}{\Omega^3(\tau)}\dfrac{\partial \hat{\hat{\widetilde H}}_1(I_*,\ph_a,\tau_*)}{\partial I}\Bigg|_{t_-}^{t_l} \\
&&{}-\eps^3\frac{3\Omega'^2(\tau)}{\Omega^5(\tau)}\dfrac{\partial \hat{\hat{\hat{\widetilde H}}}_1(I_*,\ph_a,\tau_*)}{\partial I}\Bigg|_{t_-}^{t_l} -\eps^4\frac{15\Omega'^3(\tau)}{\Omega^7(\tau)}\dfrac{\partial \hat{\hat{\hat{\hat{\widetilde H}}}}_1(I_*,\ph_a,\tau_*)}{\partial I}\Bigg|_{t_-}^{t_l} \\\\
&&{}+\eps\int\limits_{t_r}^{t_+}\frac{\partial \widetilde H_1(I_*,\ph_a,\tau_*)}{\partial I}\,\dif t -\frac{\eps}{\Omega(\tau)}\dfrac{\partial \hat{\widetilde H}_1(I_*,\ph_a,\tau_*)}{\partial I}\Bigg|_{t_r}^{t_+}-\eps^2\frac{\Omega'(\tau)}{\Omega^3(\tau)}\dfrac{\partial \hat{\hat{\widetilde H}}_1(I_*,\ph_a,\tau_*)}{\partial I}\Bigg|_{t_r}^{t_+}\\
&&{}-\eps^3\frac{3\Omega'^2(\tau)}{\Omega^5(\tau)}\dfrac{\partial \hat{\hat{\hat{\widetilde H}}}_1(I_*,\ph_a,\tau_*)}{\partial I}\Bigg|_{t_r}^{t_+} -\eps^4\frac{15\Omega'^3(\tau)}{\Omega^7(\tau)}\dfrac{\partial \hat{\hat{\hat{\hat{\widetilde H}}}}_1(I_*,\ph_a,\tau_*)}{\partial I}\Bigg|_{t_r}^{t_+}+O(\eps^{\frac32})\\
&=&\eps\int\limits_{t_-}^{t_l}\frac{\partial \widetilde H_1(I_*,\ph_a,\tau_*)}{\partial I}\,\dif t +\eps\int\limits_{t_r}^{t_+}\frac{\partial \widetilde H_1(I_*,\ph_a,\tau_*)}{\partial I}\,\dif t -\frac{\eps}{\Omega}\frac{\partial \hat{\widetilde H}_1(I_*,\ph_a,\tau_*)}{\partial I}\Bigg|_{t_-}^{t_+}+\frac{\eps}{\Omega}\frac{\partial \hat{\widetilde H}_1(I_*,\ph_a,\tau_*)}{\partial I}\Bigg|_{t_l}^{t_r}\\
&&{}+\eps^2\frac{\Omega'}{\Omega^3}\dfrac{\partial \hat{\hat{\widetilde H}}_1(I_*,\ph_a,\tau_*)}{\partial I}\Bigg|_{t_l}^{t_r}+\eps^3\frac{3\Omega'^2}{\Omega^5}\dfrac{\partial \hat{\hat{\hat{\widetilde H}}}_1(I_*,\ph_a,\tau_*)}{\partial I}\Bigg|_{t_l}^{t_r}+\eps^4\frac{15\Omega'^3}{\Omega^7}\dfrac{\partial \hat{\hat{\hat{\hat{\widetilde H}}}}_1(I_*,\ph_a,\tau_*)}{\partial I}\Bigg|_{t_l}^{t_r}+O(\eps^{\frac32}),
\end{eqnarray*}
where $\widetilde\beta$ and $\widetilde{\widetilde\beta}$ are smooth functions.

\subsubsection{Estimate of $e_7^\ph$}

Applying Lemma \ref{replaceFAR}, we easily get the estimate
$$e_7^\ph=\int\limits_{t_-}^{t_l}\sum\limits_{k=5}^{9}\eps^k\frac{\gamma^k(J,\ph,\tau)}{\omega^{2k-3}(\tau)}\,\dif t +\int\limits_{t_r}^{t_+}\sum\limits_{k=5}^{9}\eps^k\frac{\gamma^k(J,\ph,\tau)}{\omega^{2k-3}(\tau)}\,\dif t=O(\eps^{\frac32}).$$

Combining obtained estimates and the identity (\ref{phiiden}), and taking into account that $v_\pm={v_1}_\pm+O(\eps)$, we get the second formula in Theorem \ref{thm1}.

\subsubsection{Principal identity for formula (\ref{angle2})}

With results of four steps of perturbation theory, we have
\begin{eqnarray*}
\ph_*&=&\ph_l+\int\limits_{t_l}^{t_*}\omega(\tau)\,\dif t +\eps\int\limits_{t_l}^{t_*}\frac{\partial H_1(I,\ph,\tau)}{\partial I}\,\dif t\\
&=&\psi_l-\eps v_l+\int\limits_{t_l}^{t_*}\omega(\tau)\,\dif t +\eps\int\limits_{t_l}^{t_*}\frac{\partial H_1(I,\ph,\tau)}{\partial I}\,\dif t\\
&=&\psi_-+\int\limits_{t_-}^{t_l}\dot\psi\,\dif t-\eps v_l +\int\limits_{t_l}^{t_*}\omega(\tau)\,\dif t +\eps\int\limits_{t_l}^{t_*}\frac{\partial H_1(I,\ph,\tau)}{\partial I}\,\dif t\\
&=&\ph_-+\eps v_--\eps v_l+\int\limits_{t_-}^{t_*}\omega(\tau)\,\dif t +\eps\int\limits_{t_l}^{t_*}\frac{\partial H_1(I,\ph,\tau)}{\partial I}\,\dif t\\
&&{}+\int\limits_{t_-}^{t_l}\frac{\partial}{\partial J} \left(\eps\RR_1+\eps^2\RR_2+\eps^3\RR_3+\eps^4\RR_4+\eps^5\frac{\partial S_4}{\partial\tau}\right)\,\dif t+O(\eps).
\end{eqnarray*}

From this relation, similarly to the above proof, we obtain the second formula of Theorem \ref{thm2}.

\section{Numerical verification}

\subsection{Example}

We will verify previous formulas numerically using the following example suggested in \cite{1}: 
\begin{equation}\label{example_phi}
\left\{ \begin{array}{lll}
\dfrac{\dif I}{\dif t}=-\eps A(I,\tau)\cos\ph \\\\
\dfrac{\dif\ph}{\dif t}=\omega(\tau)+\eps\dfrac{\partial A(I,\tau)}{\partial I}\sin\ph
\end{array}\right.
\end{equation}

where
$$\omega(\tau)=\me^{\tau-1}-1,\quad A(I,\tau)=\frac{I\sqrt{4-I}}{\sqrt{\me^{\tau-1}+1}}.$$

We choose slow time interval $\tau\in[0,2]$, take $I(0)=1$ and consider different values $\ph(0)\in[0,2\pi]$. We have 
$$H_1(I,\ph,\tau)=A(I,\tau)\sin\ph,\quad \bar H_1=\langle H_1\rangle^\ph=0,\quad \widetilde H_1=H_1-\bar H_1=H_1,$$
$$u_1=\frac{\widetilde H_1}{\omega}=\frac{A(I,\tau)\sin\ph}{\me^{\tau-1}-1},\quad \frac{\partial u_1}{\partial I}=\frac{\partial A}{\partial I}\frac{\sin\ph}{\me^{\tau-1}-1}=-\frac{\partial v_1}{\partial\ph}, $$
$$v_1=\frac{\partial A}{\partial I}\frac{\cos\ph}{\me^{\tau-1}-1},\quad \frac{\partial A(I,\tau)}{\partial I}=\frac{\sqrt{4-I}-\dfrac{I}{2\sqrt{4-I}}}{\sqrt{\me^{\tau-1}+1}},$$
and $$I_-=1,\quad \tau_-=0,\quad \tau_*=1,\quad \tau_+=2,\quad \omega'_*=1.$$

\subsection{Theoretical value of $I_+$}

According to (\ref{action3}), we should find 
\begin{equation}\label{theorI}
I_+^{\rm{theor}}=I_-+\eps u_1(I_-,\ph_-,\tau_-)-\eps u_1(I_-,\check\ph_+,\tau_+)-\sqrt\eps\int\limits_{-\infty}^{+\infty}\frac{\partial H_1(\check I_*,\check\ph_*+\frac{\omega_*'}{2}\theta^2,\tau_*)}{\partial\ph}\,\dif\theta.
\end{equation}
We perform as follows:

1$^\circ$. Value of $\check\ph_+$.
$$\check\ph_+=\ph_-+\frac1\eps\int\limits_{\tau_-}^{\tau_+}\omega(\tau)\,\dif\tau=\ph_-+\frac1\eps\int\limits_0^2\omega(\tau)\,\dif\tau=\ph_-+\frac{\me-\me^{-1}-2}{\eps}.$$

2$^\circ$. Value of $\check{\check\ph}_*$.
$$\check{\check\ph}_*=\ph_-+\frac1\eps\int\limits_{\tau_-}^{\tau_*}\omega(\tau)\,\dif\tau=\ph_-+\frac1\eps\int\limits_0^1\omega(\tau)\,\dif\tau=\ph_--\frac1{\eps\me}.$$

3$^\circ$. Value of $\check\ph_*$.
$$\check\ph_*=\ph_-+\frac1\eps\int\limits_{\tau_-}^{\tau_*}\omega(\tau)\,\dif\tau +\frac12\sqrt\eps\int\limits_{-\infty}^{+\infty}\frac{\partial\widetilde H_1(I_-,\check{\check\ph}_*+\frac{\omega_*'}{2}\theta^2,\tau_*)}{\partial I}\,\dif\theta+\frac{\eps\ln\eps}{4\omega'_*}\frac{\partial^2\big\langle\widetilde H_1^2\big\rangle^\ph(I_-,\tau_*)}{\partial I^2}.$$
Here
\begin{eqnarray*}
&&\frac12\sqrt\eps\int\limits_{-\infty}^{+\infty}\frac{\partial\widetilde H_1(I_-,\check{\check\ph}_*+\frac{\omega_*'}{2}\theta^2,\tau_*)}{\partial I}\,\dif\theta\\
&=&\frac12\sqrt\eps\int\limits_{-\infty}^{+\infty}\frac{\ \sqrt{4-I_-}-\frac{I_-}{\ 2\sqrt{4-I_-}\ }\ }{\sqrt{\me^{\tau_*-1}+1}}\sin\big(\check{\check\ph}_*+\tfrac{\omega_*'}{2}\theta^2\big)\,\dif\theta\\
&=&\sqrt\eps\cdot\frac{\sqrt3-\frac1{2\sqrt3}}{\sqrt2}\int\limits_0^{+\infty} \left[\cos\check{\check\ph}_*\sin\tfrac{\theta^2}2 +\sin\check{\check\ph}_*\cos\tfrac{\theta^2}2\right]\,\dif\theta\\
&=&\sqrt\eps\cdot\frac5{2\sqrt6}\left(\cos\check{\check\ph}_* \int\limits_0^{+\infty} \sin\tfrac{\theta^2}2\,\dif\theta +\sin\check{\check\ph}_* \int\limits_0^{+\infty} \cos\tfrac{\theta^2}2\,\dif\theta\right)\\
&=&\sqrt\eps\cdot\frac{5\sqrt\pi}{4\sqrt6}\,(\cos\check{\check\ph}_*+\sin\check{\check\ph}_*).
\end{eqnarray*}
and
\begin{eqnarray*}
\frac{\partial^2\big\langle\widetilde H_1^2\big\rangle^\ph(I_-,\tau_*)}{\partial I^2}&=&\frac{\partial^2}{\partial I^2}\left\langle\frac{I^2(4-I)}{\me^{\tau-1}+1}\sin^2\ph\right\rangle^\ph(I_-,\tau_*)\\
&=&\frac12 \left.\frac{\partial^2}{\partial I^2}\frac{I^2(4-I)}{\me^{\tau-1}+1}\right|_{(I_-,\tau_*)}  
=\frac12 \frac{8-6I_-}{\me^{\tau_*-1}+1}
=\frac12.
\end{eqnarray*}

Therefore, $$\check\ph_*=\check{\check\ph}_*+\frac{5\sqrt{\pi\eps}}{4\sqrt6}\,(\cos\check{\check\ph}_*+\sin\check{\check\ph}_*)+\frac{\eps\ln\eps}8.$$

4$^\circ$. Value of $\check I_*$.
$$\check I_*=I_--\frac12\sqrt\eps\int\limits_{-\infty}^{+\infty}\frac{\partial H_1(I_-,\check\ph_*+\frac{\omega'_*}{2}\theta^2,\tau_*)}{\partial\ph}\,\dif\theta.$$
Here
\begin{eqnarray*}
&&\frac12\sqrt\eps\int\limits_{-\infty}^{+\infty}\frac{\partial H_1(I_-,\check\ph_*+\frac{\omega'_*}{2}\theta^2,\tau_*)}{\partial\ph}\,\dif\theta\\
&=&\frac12\sqrt\eps\int\limits_{-\infty}^{+\infty}\frac{\ I_-\sqrt{4-I_-}\ }{\ \sqrt{\me^{\tau_*-1}+1}\ }\cos\big(\check\ph_*+\tfrac{\omega_*'}{2}\theta^2\big)\,\dif\theta\\
&=&\sqrt\eps\cdot\frac{\sqrt3}{\sqrt2} \int\limits_0^{+\infty}\left[\cos\check\ph_*\cos\tfrac{\theta^2}2 -\sin\check\ph_*\sin\tfrac{\theta^2}2\right]\,\dif\theta\\
&=&\sqrt\eps\cdot\frac{\sqrt3}{\sqrt2} \left(\sqrt2\cos\check\ph_*\cdot\sqrt\frac\pi8 -\sqrt2\sin\check\ph_*\cdot\sqrt\frac\pi8\right)\\
&=&\sqrt\eps\cdot\frac{\sqrt{3\pi}}{2\sqrt2}\,(\cos\check\ph_*-\sin\check\ph_*).
\end{eqnarray*}

Therefore, $$\check I_*=1-\frac{\sqrt{3\pi\eps}}{2\sqrt2}\,(\cos\check\ph_*-\sin\check\ph_*).$$

We have also
$$\eps u_1(I_-,\ph_-,\tau_-)=\frac{\sqrt3\,\eps\sin\ph_-}{(\me^{-1}-1)\sqrt{\me^{-1}+1}},\quad \eps u_1(I_-,\check\ph_+,\tau_+) =\frac{\sqrt3\,\eps\sin\big(\ph_-+\frac{\me-\me^{-1}-2}\eps\big)}{(\me-1)\sqrt{\me+1}},$$
and
\begin{eqnarray*}
\sqrt\eps\int\limits_{-\infty}^{+\infty}\frac{\partial H_1(\check I_*,\check\ph_*+\frac{\omega'_*}{2}\theta^2,\tau_*)}{\partial\ph}\,\dif\theta &=&\sqrt\eps\int\limits_{-\infty}^{+\infty}\frac{\ \check I_*\sqrt{4-\check I_*}\ }{\ \sqrt{\me^{\tau_*-1}+1}\ }\cos\big(\check\ph_*+\tfrac{\omega_*'}{2}\theta^2\big)\,\dif\theta\\
&=&\sqrt\frac{\pi\eps}2\cdot\check I_*\sqrt{4-\check I_*} \cdot(\cos\check\ph_*-\sin\check\ph_*)
\end{eqnarray*}
Thus we have formulas for all required values for calculation of $I_+^{\rm{theor}}$  (\ref{theorI}).

\subsection{Theoretical value of $\ph_+$}

From (\ref{angle3}), taking into account that $\bar H_1=0$, $\omega'_*=1$, we get
\begin{eqnarray}
\ph_+^{\rm{theor}}&=&\ph_-+\eps v_1(I_-,\ph_-,\tau_-)-\eps v_1(I_-,\check\ph_+,\tau_+) +\frac1\eps\int\limits_{\tau_-}^{\tau_+}\omega(\tau)\,\dif\tau +\sqrt\eps\int\limits_{-\infty}^{+\infty}\frac{\partial\widetilde H_1(\check I_*,\check\ph_*+\frac{\theta^2}2,\tau_*)}{\partial I}\,\dif\theta \nonumber\\
&&{}+{\rm p.v.}\,\eps\int\limits_{\tau_-}^{\tau_+}\frac{\partial{\cal R}_2(I_-,\tau)}{\partial I}\,\dif\tau -\frac{\eps^{\frac32}\ln\eps}4\frac{\partial^3\big\langle\widetilde H_1^2\big\rangle^\ph(I_-,\tau_*)}{\partial I^3} \int\limits_{-\infty}^{+\infty}\frac{\partial H_1(I_-,\check\ph_*+\frac{\theta^2}2,\tau_*)}{\partial\ph}\,\dif\theta. \label{theorphi}
\end{eqnarray}

We perform as follows:

1$^\circ$. 
$$\ph_-+\frac1\eps\int\limits_{\tau_-}^{\tau_+}\omega(\tau)\,\dif\tau=\check\ph_+=\ph_-+\frac{\me-\me^{-1}-2}{\eps}.$$

2$^\circ$. We know that $$v_1=\frac{\partial A}{\partial I}\frac{\cos\ph}{\me^{\tau-1}-1}=\left(\sqrt{4-I}-\frac{I}{2\sqrt{4-I}}\right)\frac{\cos\ph}{(\me^{\tau-1}-1)\sqrt{\me^{\tau-1}+1}}.$$
Thus, $$\eps v_1(I_-,\ph_-,\tau_-)=\frac{\big(\sqrt3-\frac1{2\sqrt3}\big)\eps\cos\ph_-}{(\me^{-1}-1)\sqrt{\me^{-1}+1}}=\frac{\frac5{2\sqrt3}\,\eps\cos\ph_-}{(\me^{-1}-1)\sqrt{\me^{-1}+1}},$$
$$\eps v_1(I_-,\check\ph_+,\tau_+)=\frac{\big(\sqrt3-\frac1{2\sqrt3}\big)\eps\cos\check\ph_+}{(\me-1)\sqrt{\me+1}}=\frac{\frac5{2\sqrt3}\,\eps\cos\check\ph_+}{(\me-1)\sqrt{\me+1}},\qquad$$

3$^\circ$. 
\begin{eqnarray*}
&\phantom{*}&\sqrt\eps\int\limits_{-\infty}^{+\infty}\frac{\partial\widetilde H_1(\check I_*,\check\ph_*+\frac{\theta^2}2,\tau_*)}{\partial I}\,\dif\theta=\sqrt\eps\int\limits_{-\infty}^{+\infty}\frac{\sqrt{4-\check I_*}-\frac{\check I_*}{2\sqrt{4-\check I_*}}}{\sqrt{\me^{\tau_*-1}+1}}\sin(\check\ph_*+\frac{\theta^2}2)\,\dif\theta\\
&=&\sqrt\eps\cdot\frac{8-3\check I_*}{2\sqrt2\sqrt{4-\check I_*}}\cdot2\sqrt2\cdot\sqrt{\frac\pi8}(\cos\check\ph_*+\sin\check\ph_*)
=\frac{\sqrt{\pi\eps\,}(8-3\check I_*)}{\sqrt8\sqrt{4-\check I_*}}(\cos\check\ph_*+\sin\check\ph_*).
\end{eqnarray*}
Values $\check I_*, \check\ph_*$ were calculated in the previous subsection.

4$^\circ$. We have $\RR_2=-\dfrac1{2\omega}\dfrac\partial{\partial I}\big\langle\widetilde H_1^2\big\rangle^\ph=-\dfrac1{2\omega}\dfrac\partial{\partial I}\left\langle\dfrac{I^2(4-I)}{\me^{\tau-1}+1}\sin^2\ph\right\rangle^\ph$. Therefore
$$\frac{\partial \RR_2}{\partial I}=-\frac1{2\omega}\frac{\partial^2}{\partial I^2}\big\langle\widetilde H_1^2\big\rangle^\ph=-\frac1{2\omega}\cdot\frac{8-6I}{\me^{\tau-1}+1}\cdot\frac12=-\frac1{2\omega}\cdot\frac{4-3I}{\me^{\tau-1}+1}.$$
Hence,
\begin{eqnarray*}
{\rm p.v.}\,\eps\int\limits_{\tau_-}^{\tau_+}\frac{\partial{\cal R}_2(I_-,\tau)}{\partial I}\,\dif\tau &=& {\rm p.v.}\,\eps\int\limits_0^2\left(-\frac1{2\omega}\right)\frac{4-3I_-}{\me^{\tau-1}+1}\,\dif\tau \\
&=&{\rm p.v.}\,-\frac\eps2\int\limits_0^2\frac1{\me^{2\tau-2}-1}\,\dif\tau\ =\ {\rm p.v.}\,-\frac\eps2\int\limits_0^2\left(\frac{\me^{2\tau}}{\me^{2\tau}-\me^2}-1\right)\,\dif\tau\\
&=&-\frac\eps4\lim\limits_{\delta\to0+}\left(\ln|\me^{2\tau}-\me^2|\Big|_0^{1-\delta}+\ln|\me^{2\tau}-\me^2|\Big|_{1+\delta}^2\right)+\frac\eps2\tau\Big|_0^2\\\\
&=&-\frac\eps4\lim\limits_{\delta\to0+}\ln\left|\frac{\me^{2-2\delta}-\me^2}{\me^{2+2\delta}-\me^2}\right|-\frac\eps4\lim\limits_{\delta\to0+}\ln\left|\frac{\me^4-\me^2}{1-\me^2}\right|+\eps\\\\
&=&0-\frac\eps2+\eps\ =\ \frac\eps2.
\end{eqnarray*}

5$^\circ$. $$\frac{\partial^3\big\langle\widetilde H_1^2\big\rangle^\ph}{\partial I^3}=\frac{\partial^3}{\partial I^3}\left\langle\frac{I^2(4-I)}{\me^{\tau-1}+1}\sin^2\ph\right\rangle^\ph=\frac{\partial^3}{\partial I^3}\frac{I^2(4-I)}{\me^{\tau-1}+1}\cdot\frac12=\frac{-3}{\me^{\tau-1}+1}.$$
Thus
\begin{eqnarray*}
&&\frac{\eps^{\frac32}\ln\eps}4\frac{\partial^3\big\langle\widetilde H_1^2\big\rangle^\ph(I_-,\tau_*)}{\partial I^3} \int\limits_{-\infty}^{+\infty}\frac{\partial H_1(I_-,\check\ph_*+\frac{\theta^2}2,\tau_*)}{\partial\ph}\,\dif\theta\\
&=&\frac{\eps^\frac32\ln\eps}4\cdot\left(-\frac32\right)\int\limits_{-\infty}^{+\infty}\frac{\sqrt3}{\sqrt2}\cos(\check\ph_*+\frac{\theta^2}2)\,\dif\theta\\
&=&\frac{\eps^\frac32\ln\eps}4\cdot\left(-\frac32\right)\frac{\sqrt3}{\sqrt2}\cdot2\sqrt2\cdot\sqrt{\frac\pi8}(\cos\check\ph_*-\sin\check\ph_*)\\\\
&=&-\frac{3\sqrt{3\pi}}{8\sqrt2}\eps^\frac32\ln\eps\,(\cos\check\ph_*-\sin\check\ph_*).
\end{eqnarray*}

Therefore,  from (\ref{theorphi}), we obtain the estimated value of $\ph_+$:
\begin{eqnarray}
\ph_+^{\rm{theor}}&=&\ph_-+\frac{\me-\me^{-1}-2}{\eps} +\frac{\frac5{2\sqrt3}\,\eps\cos\ph_-}{(\me^{-1}-1)\sqrt{\me^{-1}+1}}-\frac{\frac5{2\sqrt3}\,\eps\cos\check\ph_+}{(\me-1)\sqrt{\me+1}}\nonumber\\
&&{}+\frac{\sqrt{\pi\eps\,}(8-3\check I_*)}{\sqrt8\sqrt{4-\check I_*}}(\cos\check\ph_*+\sin\check\ph_*)+\frac\eps2+\frac{3\sqrt{3\pi}}{8\sqrt2}\eps^\frac32\ln\eps\,(\cos\check\ph_*-\sin\check\ph_*).\qquad\label{theorphi2}
\end{eqnarray}

\subsection{Results of numerical simulation}

\begin{figure}[H]
\centering
\includegraphics[width=16cm]{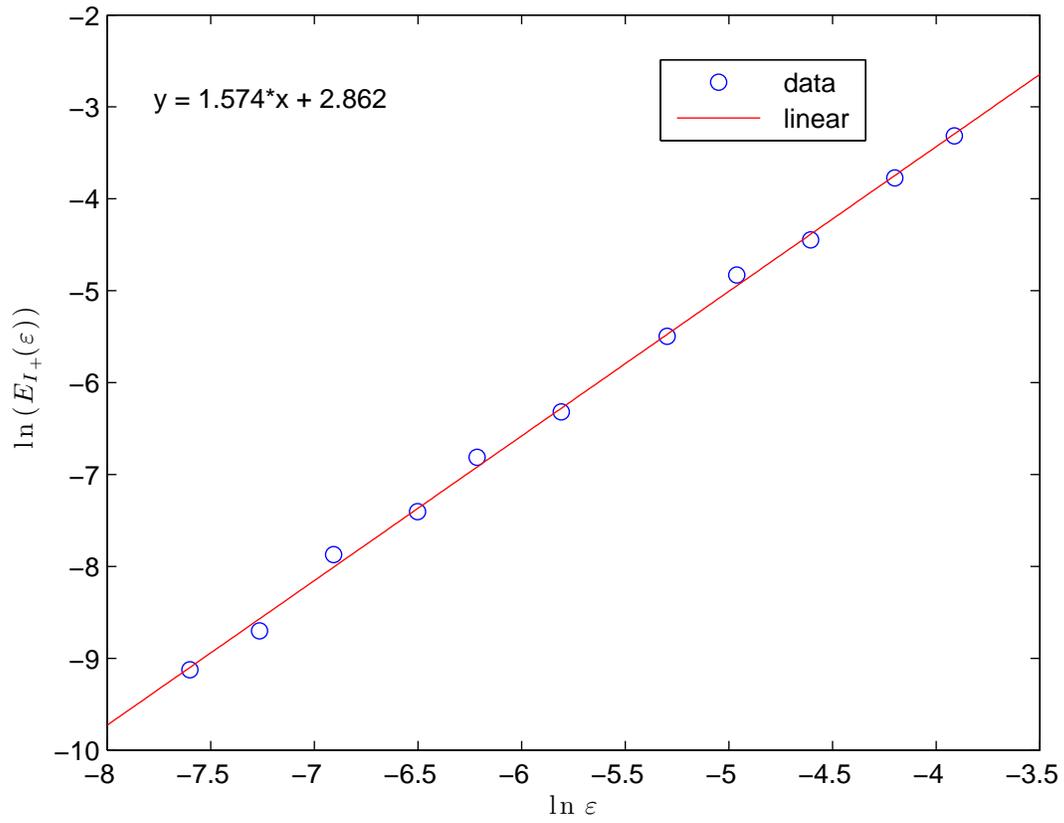}
\caption{Results of numerical simulation for action variable $I$ }
\label{figureI}
\end{figure}

\begin{figure}[H]
\centering
\includegraphics[width=16cm]{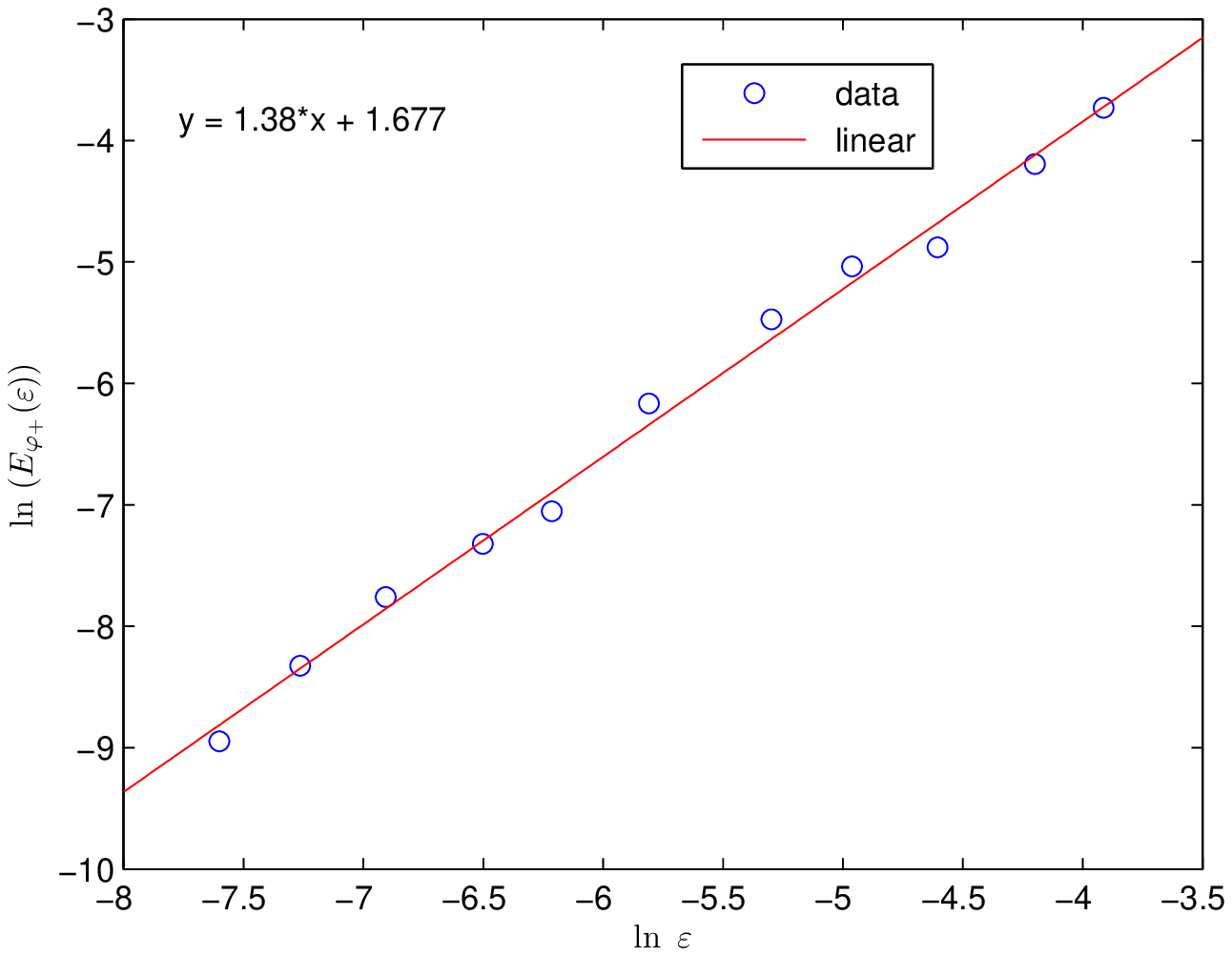}
\caption{Results of numerical simulation for angle variable $\ph$ }
\label{figurephi}
\end{figure}

Our goal here is to check numerically, if the obtained values $\ph_+^{\rm{theor}}, I_+^{\rm{theor}}$ indeed approximate actual values of $I, \ph$ on solutions of system (\ref{example_phi}) with the accuracy $O(\eps^\frac32) $, as it is predicted by    Corollary \ref{cor1}. To this end we introduce the variable 
$$
\chi = \ph-\int\limits_0^\tau\frac{\omega(\tau_1)}\eps\,\dif\tau_1=\ph-\frac1\eps\left(\me^{\tau-1}-\me^{-1}-\tau\right)
$$
and integrate numerically by 4th order Runge-Kutta algorithm the system
\begin{equation}\label{example_chi}
\left\{ \begin{array}{lll}
\dot I=-\eps A(I,\tau)\cos\Big(\chi+\dfrac{\me^{\tau-1}-\me^{-1}-\tau}\eps\Big), \\\\
\dot\chi=\eps\dfrac{\partial A(I,\tau)}{\partial I}\sin\Big(\chi+\dfrac{\me^{\tau-1}-\me^{-1}-\tau}\eps\Big),
\end{array}\right.
\end{equation}
on the slow time interval $0\le\tau\le2$. At $\tau=\tau_+= 2$ we have values  $I_+^{\rm{numer}}, \chi_+^{\rm{numer}}$ and $\ph_+^{\rm{numer}}=\chi_+^{\rm{numer}}+({\me-\me^{-1}-2})/\eps$.  These calculations are performed for 48 values of $\ph_-$ equally spaced on $[0,2\pi]$, and 11 values of $\eps$, \{0.02, 0.015, 0.01, 0.007, 0.005, 0.003, 0.002, 0.0015, 0.001, 0.0007, 0.0005\}. We calculate 
$$E_{I_+}(\eps)=\max_{\ph_-\in[0,2\pi]}|I_+^{\rm{numer}}-I_+^{\rm{theor}}|,\quad E_{\ph_+}(\eps)=\max_{\ph_-\in[0,2\pi]}|\ph_+^{\rm{numer}}-\ph_+^{\rm{theor}}|\,,$$
and plot values $\ln(E_{I_+})$ and $\ln(E_{\ph_+})$ as functions of $\ln\eps$ in Figures \ref{figureI} and \ref{figurephi}.
Linear least squares fit of the data in Figures \ref{figureI} and \ref{figurephi} gives slopes $\alpha=1.574 $ and $\beta=1.38$, respectively. The ideal results for accuracy  $O(\eps^\frac32)$ would be $\alpha=3/2$ and $\beta=3/2$. Thus the numerical simulation indicates that the accuracy is $O(\eps^\frac32)$, as expected.

\bigskip
\begin{thebibliography} {[99]}
\footnotesize

\bibitem{2} Kevorkian, J. and Cole, J. D.: \emph{Multiple Scale and Singular Perturbation Methods}. Applied Mathematical Sciences, Vol. 114. Springer-Verlag, 1996.
\bibitem{3} Chirikov, B. V.: \emph{The passage of a nonlinear oscillatory system through resonance}. Sov. Phys., Dokl. 4, 1959, pp. 390-394.
\bibitem{1} Bosley, D. L.: \emph{An improved matching procedure for transient resonance layers in weakly nonlinear oscillatory systems}. Siam J. Appl. Math, Vol. 56, No. 2, 1996, pp. 420-445.
\bibitem{4} Alekseev, P. A.: \emph{On change of action at passage through a resonance in a quasilinear Hamiltonian system}. M.Sc. thesis, Moscow University, 2007.

\end{thebibliography}
\end{document}